\documentclass[11pt]{amsart}
\usepackage{amssymb}
\usepackage{amsfonts}
\usepackage{amscd}
\usepackage{amsmath}
\usepackage{mathrsfs}
\usepackage{curves}
\usepackage{epsfig}
\usepackage[pass]{geometry}
\usepackage{graphicx}
\usepackage{verbatim}
\usepackage{latexsym}
\usepackage{eucal}

\numberwithin{equation}{section}

\newtheorem{theorem}{Theorem}[section]
\newtheorem{lemma}[theorem]{Lemma}
\newtheorem{prop}[theorem]{Proposition}
\newtheorem{cor}[theorem]{Corollary}

\def \mca {{\mathscr A}}
\def \mcb {{\mathscr B}}

\def \mcd {{\mathscr D}}
\def \mce {{\mathscr E}}

\def \mci {{\mathscr I}}

\def \mcl {{\mathscr L}}

\def \mcn {{\mathscr N}}

\def \mcr {{\mathscr R}}
\def \mcs {{\mathscr S}}

\def \mcw {{\mathscr W}}

\def \mcx {{\mathscr X}}

\def \mbr {{\mathbb R}}
\def \mbs {{\mathbb S}}

\def \loc {\text{loc}}

\def \re {\operatorname{Re}}

\def \defeq {\stackrel{\operatorname{def}}{=}}
\def \beqq {\begin{equation}}
\def \eeqq {\end{equation}}
\def \WF {\text{WF}}

\def \bpf {\begin{proof}}
\def \epf {\end{proof}}
\def \beq {\begin{equation*}}
\def \eeq {\end{equation*}}

\def \eps {\epsilon}   
\def \la {\lambda}   
\def \La {\Lambda}

\def \p {\partial}
\def \ha {\frac{1}{2}}

\def \tilde {\widetilde}

\begin{document}
\title[Metal artifacts in X-ray tomography]{Quantitative analysis of metal artifacts in X-ray tomography }
\author{Benjamin Palacios	}
\address{Benjamin Palacios
\newline
\indent Department of Mathematics, University of Washington}
\email{bpalacio@uw.edu}

\author{Gunther Uhlmann} 
\address{Gunther Uhlmann
\newline
\indent Department of Mathematics, University of Washington,
\newline
\indent and Institute for Advanced Study, the Hong Kong University of Science and Technology}
\email{gunther@math.washington.edu}

\author{Yiran Wang}
\address{Yiran Wang
\newline
\indent Department of Mathematics, University of Washington}
\email{wangy257@math.washington.edu}

\begin{abstract} 
In X-ray CT scan with metallic objects, it is known that direct application of the filtered back-projection (FBP) formula leads to streaking artifacts in the reconstruction. These are characterized mathematically in terms of wave front sets in \cite{Seo}. In this work, we give a quantitative microlocal analysis of such artifacts. We consider metal regions with strictly convex smooth boundaries and show that the streaking artifacts are conormal distributions to straight lines tangential to at least two boundary curves. For metal regions with piecewise smooth boundaries, we analyze the streaking artifacts especially due to the corner points. Finally, we study the reduction of the artifacts using appropriate filters. 
\end{abstract}

\maketitle

\section{Introduction}
\begin{figure}[htbp]
\centering
\includegraphics[scale =.6 ]{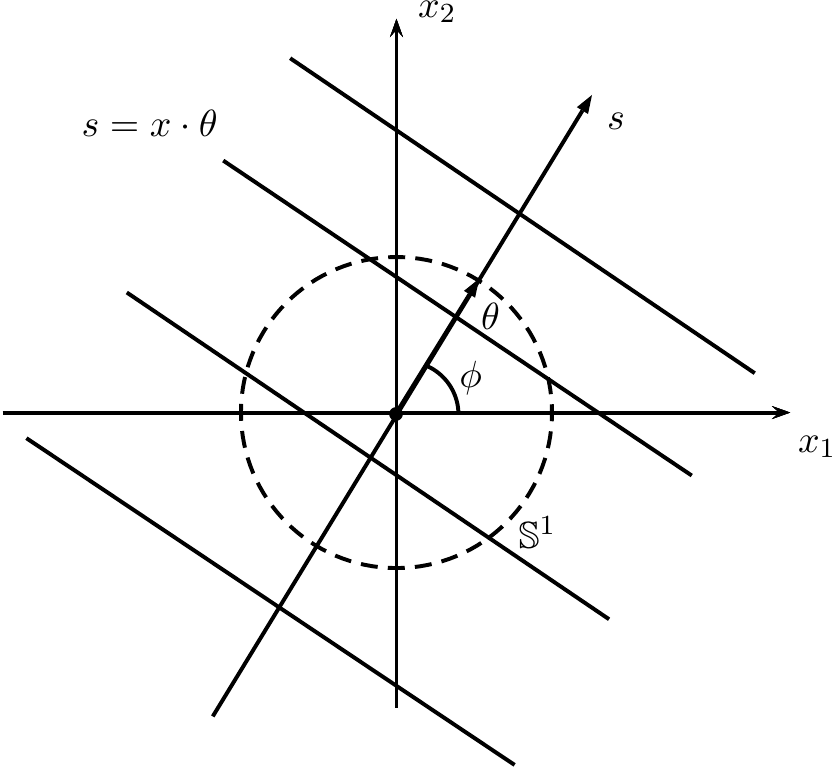}
\caption{Radon transform on $\mbr^2.$ The transform is to integrate $f(x)$ along the lines $s = x\cdot \theta.$}
\label{figradon}
\end{figure}

X-ray computed tomography (CT) is widely used in medical and dental imaging. In CT scan, X-ray projection data $P(s, \phi)$, also known as the sinogram, is collected after X-ray beams passing through in different directions, see Figure \ref{figradon}. More precisely, let $f(x)$ denote the attenuation coefficients of X-rays corresponding the object being imaged. For most human tissues, $P(s, \phi) = \mcr f$, where $\mcr $ is the X-ray transform on $\mbr^2$ (or the Radon transform):  
\beqq\label{eqradon}
\begin{gathered}
\mcr f(s, \phi) = \int_{\mbr^2} \delta(  x \cdot \theta -s) f(x) d  x, \\
  x = (x_1, x_2) \in \mbr^2, \ \ \theta = (\cos \phi, \sin \phi), \ \ \phi \in (-\pi, \pi], \ \ s\in \mbr.
\end{gathered}
\eeqq
Here, $(s, \theta)$ are local coordinates for the cylinder $M \defeq \mbr \times \mbs^1$. We also use $(s, \phi)$ for the coordinates when convenient. The (inverse) problem is to find $f(x)$ from $P(s, \theta)$. For Radon transform, we have the well-known filtered back projection (FBP) formula, see e.g.\ \cite{Na}
\beq
f = \frac{1}{4\pi} \mcr^* \mci^{-1} \mcr f, \ \ f\in \mce'(\mbr^2),
\eeq
where 
\beq
\begin{gathered}
\mcr^*h(  x) = \int_{-\pi}^\pi h(\phi,   x\cdot \theta) d\phi \text{ is the adjoint of $\mcr,$}\\
\text{ and } \mci^{-1}(g)(s) = \frac{1}{2\pi} \int_\mbr\int_\mbr e^{i(s-s') w} g(s')|w|ds' dw \text{ is the Riesz potential. }
\end{gathered}
\eeq
In practice, this algorithm works well if the CT scanned data $P(s, \phi)$ belongs to the range of $\mcr$ with domain $\mce'(\mbr^2)$, which is approximately true for most human tissues. 

In CT-scan containing metallic objects, the data $P(s, \phi)$ may not belong to the range of the Radon transform and direct application of the FBP formula leads to streaking artifacts. The artifacts appear as line segments in the reconstructed image of $f$. The understanding of the cause of the artifacts lead to challenging problems and there is a huge literature on the reduction methods. We refer the reader to \cite[Section 5]{Seo1}, \cite{BK} for an overview and \cite{ZDL} for some recent progress. However, the main focus of this work is the  mathematical study of the artifacts from a quantitative point of view. 

Several causes have been identified to account for the mismatch. The cause we consider in this work is that the attenuation of metallic objects varies largely with respect to the energy level $E$ of the X-ray, but we remark that our analysis applies to other causes such as the crompton scattering. We use the model derived in the very nice paper \cite{Seo} which we summarize in what follows. Let $D\subset \mbr^2$ denote the region of metallic objects and $\chi_D$ be the characteristic function of $D$. Suppose that $E$ is in an energy window $[E_0-\eps, E_0+\eps]$ for $\eps>0$ small. The attenuation coefficient $f_E$ of the object being imaged at energy level $E$ can be described as
\beq
f_E = f_{E_0} + \alpha (E-E_0) \chi_D,
\eeq
where $\alpha$ is a constant. As in \cite{Seo}, we make the technical assumptions that (i) $f_E\in \mce'(\mbr^2)$ and (ii) $f_{E_0}(x) \geq C \sup_{y\in \bar D^c} f_{E_0}(y)$ for $x\in D$ with some $C> 1$. Here, one should think of $\alpha$ as an approximation of the derivative $\frac{\p f_E}{\p E}$ in $D$, see \cite[equation (2.4)--(2.7)]{Seo}. For normal tissues, one has $\frac{\p f_E}{\p E} \approx 0$ so that $f_E\approx f_{E_0}.$ However, for metallic objects, the derivative is not small. 

Now let $\eta(E)$ denote the fractional energy at photon energy $E$. The X-ray data is given by
\beq
\begin{split}
P(s, \phi) &= -\ln [\int_{E_0-\eps}^{E_0+\eps} \eta(E) \exp\{-\mcr f_E(\phi, s)\} dE] \ \ \text{ (Beer's law)}\\
& \approx \mcr f_{E_0} - \ln \big(\frac{\sinh (\alpha \eps \mcr \chi_D)}{\alpha \eps \mcr \chi_D}\big),
\end{split}
\eeq
where the second line is obtained by taking an approximation $\eta(E) = 1/2\eps$, see \cite[equation (2.11)]{Seo}. Let $P_{MA} = P(s, \phi) - \mcr f_{E_0}$ be the mismatch. Throughout the rest of the paper, we shall work with an approximation of $P_{MA}$
\beq
P_{MA, N} = \sum_{k = 1}^N \frac{(-1)^k}{k} [\sum_{n = 1}^N \frac{(\alpha \eps)^{2n}}{(2n+1)!}(\mcr \chi_D)^{2n}]^k, \ \ N \geq 1. 
\eeq
It is shown in \cite[Proposition 2.2]{Seo} that $P_{MA, N}$ converges to $P_{MA}$ in $H^t(M)$ for $t\in (0, \ha)$. For convenience, we shall abuse the notations $P_{MA, N}$ and $P_{MA}$. For this data, direct application of the  FBP formula gives the reconstruction formula  
\beqq\label{eqfct}
\begin{gathered}
f_{CT}(  x) = \frac{1}{4\pi} \mcr^* \mci^{-1} P = f_{E_0}(x) + f_{MA}(x), \ \ \text{where }\\
f_{MA}(  x) = \frac{1}{4\pi} \mcr^* \mci^{-1}[\sum_{k = 1}^N \frac{(-1)^k}{k} [\sum_{n = 1}^N \frac{(\alpha \eps)^{2n}}{(2n+1)!}(\mcr \chi_D)^{2n}]^k]. 
\end{gathered}
\eeqq
To our knowledge, the first microlocal description of the metal artifacts was introduced in \cite{Seo} from a qualitative perspective. 
In \cite[Section 3]{Seo}, the authors defined the streaking artifacts using the notion of wave front set and identified that $f_{MA}$ is the distribution containing the streaking artifacts. Moreover, the authors obtained conditions under which the artifacts would appear. Based on these results, numerical methods are proposed in \cite{Seo1} to reduce such artifacts. Our goal in this work is to give a precise quantitative description of the streaking artifacts using microlocal methods by determining their strength in terms of their order as singularities in the reconstruction. We should warn the reader that in this work, the word {\em streaking artifact} is referred to both the line segment in the reconstructed image of $f_{CT}$ which is an geometric object, and the distribution in $f_{CT}$ associated with such artifacts which is a "function". However, it should be clear from the context which one we are referring to.  

It is obvious that $f_{MA}$ is a nonlinear function of $\mcr\chi_D$, however we must emphasize that it is indeed the nonlinear interactions of the singularities in $\mcr\chi_D$ that produces the streaking artifacts. The geometry of the metal regions plays an important role. For metal regions with smooth boundaries, the streaking artifacts intersect  the boundary of metallic objects in a particular way so that their singularities can be described using the notion of paired Lagrangian distributions, see Theorem \ref{main0} and \ref{main1} for the detail. We characterize the streaking artifacts in both the data (sinogram) and the reconstruction. More importantly, away from the boundaries of $D, $ the streaking artifacts are conormal distributions and we determine their strength (order). Our analysis leads us to construct appropriate filters which makes the artifacts smoother, see Theorem \ref{main3}. Furthermore, our analysis can be applied to metal regions with piecewise smooth boundaries, see Theorem \ref{main2}. This is mentioned in \cite{Seo} but not addressed with much detail. The situation is interesting but more complicated, because streaking artifacts can also be generated from the corner points. 

We remark that characterizing streaking artifacts using wave front sets has been considered for quantitative susceptibility mapping (QSM) in \cite{Seo0}. Microlocal techniques are used to reduce the effects in \cite{PUW}. Although we do not pursue it here, we point out that it is worth studying the generalization of the problem to geodesic ray transforms on Riemannian manifolds, where the artifacts should follow geodesic rays. This arises in ultrasound. 

The paper is organized as follows. We start with some preliminary analysis in Section \ref{pre}. In Section \ref{secnon} and \ref{secfull}, we consider metal regions with smooth boundaries. We start with analyzing the singularities in the first term of $f_{MA}$ in Section \ref{secnon} to show the nonlinear effects. Then we study the singularities of the full $f_{MA}$ in Section \ref{secfull}. We deal with metal regions with piecewise smooth boundaries in Section \ref{secpiece}. Finally, we consider the reduction of the streaking artifacts using appropriate filters. In Appendix \ref{secapp}, we recall the definition and basics of conormal and paired Lagrangian distributions for readers' convenience.

\section*{Acknowledgement}
GU was partly supported by NSF, a Walker Family Endowed Professorship at UW and a Si-Yuan Professorship at HKUST. BP was partially supported by NSF grant  DMS-1265958.

\section{Preliminaries}\label{pre}

We assume that the metal region $D = \bigcup_{j=1}^J D_j$, $J \geq 1$ where $D_j$ are simply connected, pair-wisely disjoint bounded domains in $\mbr^2$ with boundary $\Sigma_j \defeq \p D_j$. We first work with a simpler geometrical setup by assuming that 
\beqq\label{eqap1}
\begin{split}
\text{ $\Sigma_j$ are strictly convex smooth curves.}
\end{split}\tag{A1}
\eeqq
Recall that a smooth curve $\gamma \subset \mbr^2$ is {\em strictly convex} if any straight lines intersect $\gamma$ at most at two points. We make few remarks. 
\begin{enumerate}
\item We assumed $D_j$ to be simply connected for simplicity. Our analysis applies to the case when $D_j$ are connected and the boundary $\Sigma_j$ is the union of finitely many connected $\Sigma_j^i, i = 1, \cdots, N_j$ satisfying (A1). 
\item  $\Sigma_j$ will be generalized to piecewise smooth curves in Section \ref{secpiece}. 
\item  (A1) implies that for distinct $\Sigma_j$ and $\Sigma_k$, there are only finitely many lines tangent to both of them.  It excludes the possibility that there are infinitely many tangent lines and they  converge, in which case the resulting singularities are expected to be complicated.  
\end{enumerate}

We denote $\Sigma = \bigcup_{j = 1}^J \Sigma_j.$ The characteristic function $\chi_{D_j}, j = 1, \cdots, J$ have Heaviside type singularities at $\Sigma_j$ and we describe them using H\"ormander's notion of Lagrangian distributions, see Appendix \ref{secapp} for a brief recall of the notation and basics of such distributions. We have $\chi_{D_j} \in I^{-1}(N^*\Sigma_j)$. (The order of the symbol of $\chi_{D_j}$ is $-1 = \mu + \frac n4 - \frac k2$ with $\mu$ the order of the Lagrangian distribution, $n = 2$ the dimension of the ambient space and $k = 1$ the co-dimension of $\Sigma_j$ in $\mbr^2$. So we have $\mu=-1.$)  Finally, we have 
\beq
\chi_{D} = \sum_{j = 1}^J \chi_{D_j}  \in \sum_{j = 1}^J I^{-1}(N^*\Sigma_j)  = I^{-1}(N^*\Sigma).
\eeq

Let us consider $\mcr \chi_{D_j}, j = 1, \cdots J$. It is known that $\mcr: \mce'(\mbr^2)\rightarrow \mcd'(M)$ is an elliptic Fourier integral operator. As $\chi_{D_j}$ is a Lagrangian distribution, we shall apply H\"ormander's FIO theory to analyze their composition and show that the resulting distribution $\mcr\chi_{D_j}$ is still conormal. Using local coordinates $(s, \theta)$ and $x$ in \eqref{eqradon}, we write the Schwartz kernel of $\mcr$, denoted by $K_\mcr$, as an oscillatory integral
\beq
K_\mcr(s, \theta, x) = \frac{1}{(2\pi)^\ha} \int_\mbr  e^{i(  x \cdot \theta -s)\la} d\la.
\eeq
The phase function is $\phi(s, \theta, x; \la) = (x \cdot \theta -s)\la$ so the associated Lagrangian submanifold is 
\beq
\begin{split}
\La &= \{(s, \theta,  x; d_{s, \theta}\phi, -d_{  x} \phi) \in T^*(M \times \mbr^2)\backslash 0 : d_\la \phi = 0\}\\
 &= \{(s, \theta, x; -\la, \la   x, -\la \theta): s = x\cdot \theta,  \la \in \mbr,   x\in \mbr^2, \theta\in \mbs^1\}.
 \end{split}
\eeq
Therefore $K_\mcr \in I^{-\ha}(\La)$. We denote the homogeneous canonical relation by 
\beqq\label{eqcaC}
C \defeq \La' = \{(s, \theta, -\la, \la   x; x,  \la \theta): s = x\cdot \theta,  \la \in \mbr,   x\in \mbr^2, \theta\in \mbs^1\} \subset T^*M \times T^*\mbr^2.
\eeqq
In this case, $\La$ (or $C$) gives us information about how the operator $\mcr$ moves singularities from the phase space in $\mbr^2$, to the phase space of $M$. 
We know that $\chi_{D_j}$ is a Lagrangian distribution associated with  $ N^*\Sigma_j$. 
Let $C_j \defeq (N^*\Sigma_j)'\subset T^*\mbr^2$ be the canonical relation. The two homogeneous canonical relations $C, C_j$ intersect transversally and their composition is 
\beq
 C\circ C_j = \{(s, \theta; -\la, \la   x) \in T^*M\backslash 0 :   s = x\cdot \theta, \ \ x\in \Sigma_j, \la \in \mbr\backslash 0, \theta \in N_x^*\Sigma_j \cap \mbs^1\}.
\eeq
This is a Lagrangian submanifold of $T^*M$, but we claim that it is a conormal bundle under our assumptions. In fact, the projection of $C\circ C_j$ to $M$ is injective by assumption (A1) and the projection is 
\beqq\label{defsji}
S_j \defeq \{(s, \theta) \in M : s = x\cdot \theta, \ \ \theta\in N_x^*\Sigma_j \cap \mbs^1, x\in \Sigma_j \}.
\eeqq
These are co-dimension one submanifolds of $M$ and $(s, \theta)\in S_j$ implies $(-s, -\theta)\in S_j$. We see that $C\circ C_j =  N^*S_j$. The fact that $C\circ C_j$ are conormal is what allows us to do the analysis of the product of such distributions. One can apply H\"ormander's clean FIO composition theorem \cite[Theorem 25.2.3]{Ho4} to conclude that 
\begin{lemma}
Under assumption (A1) and with $S_j$ defined in \eqref{defsji}, $\mcr \chi_{D_j} \in I^{-\frac 32}(N^*S_j)$ is a conormal distribution. We denote $N^*S  \defeq \bigcup_{j=1}^J  N^*S_j$ and conclude that $\mcr \chi_{D} \in I^{-\frac 32}(N^*S)$. 
\end{lemma}
Roughly speaking, the singularities in the sinogram of $\chi_D$ are contained in the conormal bundle of the curve $S$. Moreover, the strength of the singularities is reduced from order $-1$ to $-3/2$.

\section{Microlocal analysis of the nonlinear effects}\label{secnon}
We start by analyzing the singularities in the first term of $f_{MA}$, that is 
\beqq\label{eqfma}
f_{MA, 1}(x)  \defeq  -\frac{1}{4\pi} \mcr^* \mci^{-1}  \frac{(\alpha \eps)^{2 }}{3!}(\mcr \chi_D)^{2 }.
\eeqq 
 We shall give a clear description of the singularities in $f_{MA, 1}$ using the notion of conormal and paired Lagrangian distributions. We address this term separately because it already demonstrates the nonlinear effects and reveals the main feature of the singularities of the artifacts.
 
Consider the nonlinear term in \eqref{eqfma} 
\beq
[\mcr(\chi_D)]^2 = \sum_{i = 1}^J [\mcr(\chi_{D_i})]^2 + 2\sum_{1\leq i< j\leq J} \mcr(\chi_{D_i})\mcr(\chi_{D_j}), \ \ J\geq 2.
\eeq
If $J = 1$, we would only have the first summation term. Because $\mcr \chi_{D_j} \in I^{-\frac 32}(N^*S_j)$, we conclude from Corollary \ref{lmprod} later in Section \ref{secfull} that $[\mcr \chi_{D_i} ]^2 \in I^{-\frac 32}(N^*S_j)$. Therefore $\WF([\mcr \chi_{D_i} ]^2) \subset N^*S_j$ does not produce new singularities. However, if the singular support of $\mcr\chi_{D_i}, \mcr\chi_{D_j}$, denoted by $S_i, S_j$ respectively, intersect transversally at $S_{ij}$, then $\mcr \chi_{D_i} \mcr \chi_{D_j} $ has new singularities at $S_{ij}$ by a wave front analysis. The resulting singularities can be described precisely using the notion of {\em paired Lagrangian distributions}. We refer the reader to Appendix \ref{secapp} for the definition and basic properties of such distributions.

\begin{figure}[htbp]
\centering

\includegraphics[scale=.7]{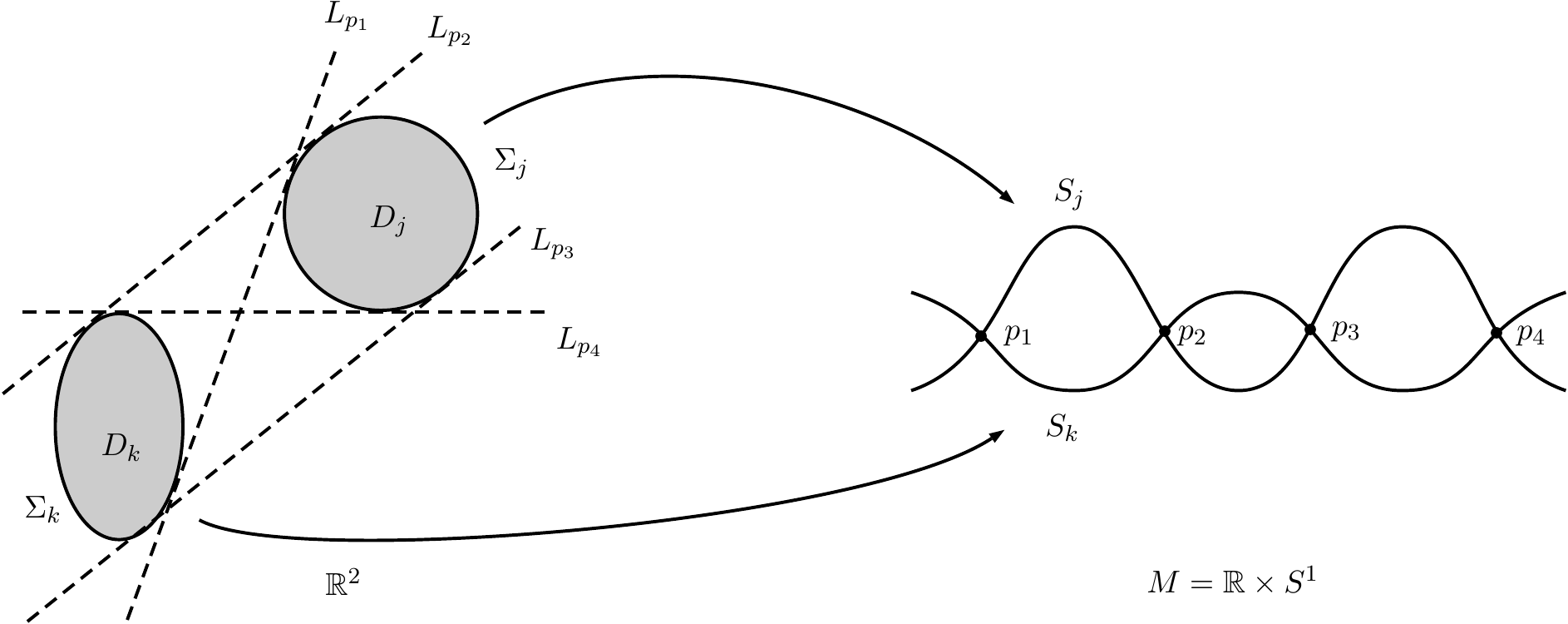}
\caption{A schematic view of singularities of $\chi_{D_j}, \chi_{D_k}$ on $\mbr^2$ (Left figure) and those of $\mcr\chi_{D_j}, \mcr \chi_{D_k}$ on $M$ the sinogram (Right figure). The singular supports $\Sigma_j, \Sigma_k$ are mapped to two curves $S_j, S_k$ respectively. Each intersection point $p_i$ corresponds to a tangent line (the dashed lines) $L_{p_i}$ representing the streaking artifacts.}
\label{figsino}
\end{figure}

Consider the intersections of $S_j, j = 1, \cdots, J$, see Figure \ref{figsino}.
\begin{lemma}\label{lmgeo}
By our assumptions on $D_j, j = 1, \cdots, J$ and their boundaries $\Sigma_j$, we have 
\begin{enumerate}
\item For $j, k = 1, \cdots, J$, $j\neq k$, $S_j$ intersect $S_k$ transversally at a finite point set $S_{jk} \subset M$; 
\item For each $p \in S_{jk}$,  there is a straight line $L_{p}\subset \mbr^2$ tangent to both $\Sigma_j$ and $\Sigma_k$. 
\end{enumerate}
\end{lemma}
\bpf
Suppose that $S_j\cap S_k \neq \emptyset$ for some $j, k = 1, \cdots, J$.  Let $p_0  = (s_0, \theta_0)$ be a point in the intersection set. By definition \eqref{defsji}, we know that there exists $x_j \in \Sigma_{j}$ and $x_k\in \Sigma_k$ so that 
\beq
s_0 = \theta_0 \cdot x_j \text{ and } s_0 = \theta_0 \cdot x_k. 
\eeq
Therefore, $x_j, x_k$ lie on the straight line 
\beqq\label{eqline}
L_{p_0} = \{x\in \mbr^2: s_0 = \theta_0 \cdot x, \ \ p_0 = (s_0, \theta_0)\in M\}.
\eeqq
It is obvious that $L_{p_0}$ is tangent to both $\Sigma_j$ and $\Sigma_k$. Because $\Sigma_j, \Sigma_k$ are strictly convex, there exits lines tangent to both of them. Thus $S_j\cap S_k\neq \emptyset.$

Next, we show that the intersection is transversal. Notice that the intersection $T_{p_0}S_{j}\cap T_{p_0}S_k$ is either (1): a one dimensional space, or (2): the zero vector.  In case (2), the intersection is transversal. Case (1) implies that the normal vectors are linearly dependent. The normal vectors to $S_j, S_k$ are spanned by $(-1, x_j), x_j\in \Sigma_j$ and  $(-1, x_k), x_k\in \Sigma_k$ respectively. If they are linearly dependent, we get $x_k = x_l$ which implies that $\Sigma_j\cap \Sigma_k\neq \emptyset$. This contradicts to the assumption that $D_j$ are pair-wisely disjoint.  Hence the intersection must be transversal. 
\epf

We introduce some notations. For $1\leq i< j\leq J$, we know that the intersection set $S_{ij}$ are point sets. We let $\mcs$ be the union of such $S_{ij}$. Corresponding to $S_{ij}$, we define $L_{ij} = \{L: \text{$L$ is a line tangent to $D_i$ and $D_j$}\}$. Finally, we let $\mcl$ to be the union of $L_{ij}$, that is 
\beqq\label{defL}
\mcl = \{L: \text{$L$ is a line tangent to $D_i$ and $D_j, 1\leq i < j \leq J$}\}.
\eeqq
The importance of $\mcl$ is that they represent the streaking artifacts as shown below. We denote 
\beq
N^*S_{ij} = \bigcup_{p\in S_{ij}} N^*p,\ \  N^*\mcs = \bigcup_{1\leq i< j\leq J} N^*S_{ij} \text{ and } N^*L_{ij} = \bigcup_{l\in L_{ij}}N^*l, \ \ N^*\mcl = \bigcup_{1\leq i< j\leq J} N^*L_{ij}.
\eeq
Notice that these notations denote the union of the conormal bundles instead of the actual conormal bundle of the union. 

Since $S_i, S_j$ intersect transversally at $S_{ij}$, using \cite[Lemma 1.1]{GrU93}, we obtain that locally near any $p\in S_{ij}$
\beqq\label{eqdij}
\mcr(\chi_{D_i})\mcr(\chi_{D_j}) \in I^{-\frac 32, -\frac 32 + \ha}(N^*p, N^*S_i) + I^{-\frac 32, -\frac 32+ \ha}(N^*p, N^*S_j).
\eeqq
In fact, it follows from the proof of the lemma that the symbol of $\mcr(\chi_{D_i})\mcr(\chi_{D_j})$ at $N^*p$ is non-vanishing. We remark that our order here is the order of Lagrangian distributions  instead of the order of symbols used in \cite[Lemma 1.1]{GrU93}. The conversion can be found in the paragraph below equation (1.4) of \cite{GrU93}. We see that the wave front set is contained in $N^*p\cup N^*S_i\cup N^*S_j$ and 
\beq
\mcr(\chi_{D_i})\mcr(\chi_{D_j}) \in I^{-\frac 52}(N^*p \backslash N^*(S_i\cup S_j))
\eeq
locally near $p$. Of course this applies to other points in $S_{ij}$. We show that this  distribution  carries the streaking artifacts.

\begin{prop}\label{propfbp}
Away from $N^*\Sigma_i\cup N^*\Sigma_j$, we have  that $\mcr^*\circ \mci^{-1} \big( \mcr(\chi_{D_i})\mcr(\chi_{D_j})\big) \in I^{-2}(N^*L_{ij})$ for $1\leq i < j \leq J$ and the principal symbol is non-vanishing. 
\end{prop}
\bpf
$\mci^{-1}$ is an elliptic pseudo-differential operator of order $1$. Thus 
\beq
\mci^{-1}\big( \mcr(\chi_{D_i})\mcr(\chi_{D_j}) \big) \in \sum_{p\in S_{ij}}\big[ I^{-\ha, -1}(N^*p, N^*S_i) + I^{-\ha, -1}(N^*p, N^*S_j) \big]
\eeq
and the principal symbol at $N^*p$ is non-vanishing. Also, we know that $\mcr^*$ is an elliptic FIO of order $-\ha$. Let $C^*$ be the canonical relation given by
\beq
C^* = \{(x, \theta, s;  \la \theta, \la x, -\la): s = x\cdot\theta, \la \in \mbr\backslash 0,   x\in \mbr^2, \theta\in \mbs^1.\}
\eeq
Then we check that 
\beq
C^*\circ N^*p = \{(x;  \la \theta_0) : \theta_0\cdot x = s_0, (\theta_0, s_0) = p,  \la \in \mbr\backslash 0\} = N^*L_p,
\eeq 
where $L_p$ is the line tangent to $\Sigma_i, \Sigma_j$ corresponding to $p$ as defined in \eqref{eqline}. We check that $C^*\circ N^*S_{i} = C^*\circ C\circ N^*\Sigma_i =  N^*\Sigma_i$.  So we get 
\beq
\mcr^*\circ \mci^{-1}\big(\mcr(\chi_{D_i})\mcr(\chi_{D_j})\big)\in \sum_{p\in S_{ij}} \big[ I^{-1, -1}(N^*L_{p}, N^*\Sigma_i) + I^{-1, -1}(N^*L_{p}, N^*\Sigma_j) \big],
\eeq
and the principal symbol at $N^*L_p$ is non-vanishing, see for example \cite{GU}. This implies the conclusion.  
\epf

The result implies that under the assumptions on the metallic region, $f_{MA, 1}$ is a paired Lagrangian distribution. According to \cite[Definition 3.2]{Seo}, the straight lines $L_{ij}$ are {\em streaking artifacts} in the sense of wave front sets. The proposition above states that they are represented by conormal distributions. We summarize the main result. 
 
\begin{theorem}\label{main0}
Consider $f_{MA, 1}$ in $\mce'(\mbr^2)$ defined in \eqref{eqfma}. Under assumptions (A1), we have that away from $\Sigma = \p D$, the streaking artifacts $f_{MA, 1} \in I^{-2}(N^*\mcl)$ with $\mcl$ defined in \eqref{defL}.  
\end{theorem}

In fact, we have a strong conclusion that if all of the $L_{ij} \in \mcl, 1\leq i < j \leq J$ are distinct which happens if there is no streak line tangent to more than two metal regions, then we have 
\beq
\WF(f_{MA, 1}) = N^*\mcl.
\eeq
However, it is not clear whether this is true in general, especially when higher order nonlinear terms are considered, see Theorem \ref{main1}. 

\begin{figure}[htbp]
\centering
\includegraphics[scale=.7]{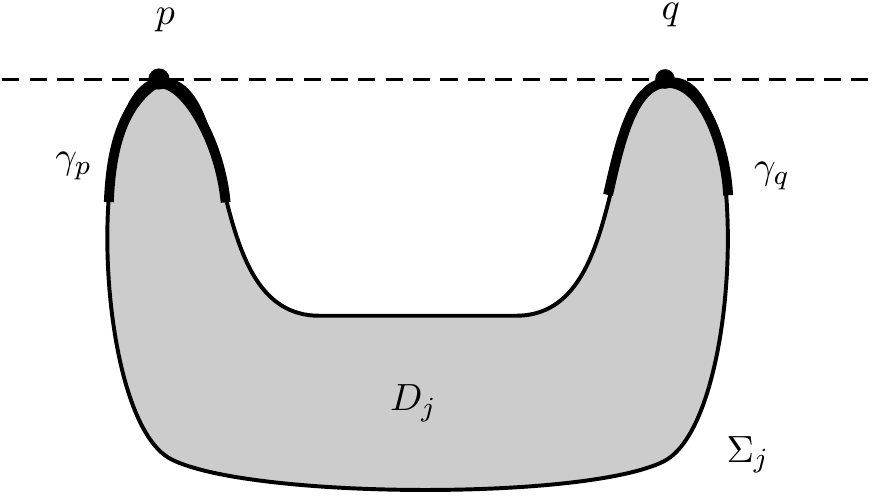}
\caption{Locally strictly convex regions and the streaking artifacts}
\label{figloccon}
\end{figure}

To conclude this section, we remark on the strictly convexity assumption in (A1). In principle, we can relax the assumption to locally strictly convexity, by which we mean that for any $p\in \Sigma_j,$ there is a neighborhood $\gamma_p$ of $p$ on $\Sigma_j$ so that $\gamma_p$ is strictly convex, see Figure \ref{figloccon}. In this case, one can introduce a partition of unity $\psi_k, k = 1, \cdots, K$ for $D_j$ so that the singular support of $\psi_k \chi_{D_j}$ is strictly convex. Since $\chi_{D_j} = \sum_{k = 1}^K \psi_k \chi_{D_j}$, our analysis goes through with some minor modifications. In particular, streaking artifacts may appear when there is a line tangent to $\gamma_p$ and $\gamma_q$ as shown in Figure \ref{figloccon}.

The situation which is still ruled out by the locally strictly convex assumption is that $\Sigma_j$ contains a line segment.  In this case, the main complexity is that $C\circ N^*\Sigma_j$ is a Lagrangian submanifold of $T^*M$ but not a conormal bundle (the projection of $C\circ N^*\Sigma_j$ to $M$ is not injective).  We hope to address this case in the future. 

\section{Metal regions with smooth boundaries}\label{secfull}
We study the singularities in the full $f_{MA}$ defined in \eqref{eqfct}. Essentially, we need to study the following terms 
\beqq\label{eqhigher}
[(\mcr \chi_D)^2]^k = [\sum_{i = 1}^J (\mcr \chi_{D_i})^2 + 2\sum_{1\leq i< j\leq J} (\mcr \chi_{D_i} \mcr \chi_{D_j})]^k, \ \ k = 2, 3, \cdots.
\eeqq
We shall prove in Corollary \ref{lmprod} that  $(\mcr \chi_{D_i})^2 \in I^\mu(N^*S_i)$ is conormal so that the term in the square bracket of \eqref{eqhigher} belongs to 
\beqq\label{defs}
\mci \defeq \sum_{i = 1}^J I^{\mu}(N^*S_i) + \sum_{S_{ij}\in \mcs} I^{\mu, \mu + \ha}(N^*S_{ij}, N^*S_i)
\eeqq
combining with our analysis in the previous section. Hereafter, for convenience we let $p = \mu, l = \mu + \ha$ and we use the notation 
\beq
I^{p, l}(N^*S_{ij}, N^*S_i) \defeq \sum_{q\in S_{ij}} [I^{p,l}(N^*q, N^*S_i) + I^{p, l}(N^*q, N^*S_j)].
\eeq
We need to analyze the product $uv$ of the following types: 
\begin{enumerate}
\item $u \in I^{p, l}(N^*S_{ij}, N^*S_i), v\in I^{\mu}(N^*S_j)$; 
\item $u \in I^{p, l}(N^*S_{ij}, N^*S_i), v\in I^{\mu}(N^*S_i)$;
\item $u \in I^{p, l}(N^*S_{ij}, N^*S_i), v\in I^{p, l}(N^*S_{ij}, N^*S_i)$;
\item $u \in I^{p, l}(N^*S_{ij}, N^*S_i), v\in I^{p, l}(N^*S_{ij}, N^*S_j)$.
\end{enumerate}
Some of these multiplications are studied in \cite{GrU93}. We will simplify the analysis by sacrificing the optimality of the results.  {\em Our goal is to show that $\mci$ is an algebra under distribution multiplication for the exponents below}
\beq
\mu \leq -\frac 32, \ \ p = \mu \leq -\frac 32, \ \  l = \mu + \ha\leq -1.
\eeq
The results we prove in the follows do not apply to general exponents and we do not pursue the optimal order for the resulting distributions $uv$. These considerations save us from some technical discussions.  

Before the proof, let us recall a useful fact of paired Lagrangian distribution that
\beqq\label{eqipl}
\bigcap_{l} I^{p, l}(\La_0, \La_1) = I^p(\La_1),
\eeqq
see \cite[Proposition 6.2]{GU}. Then  case (2) can be reduced to case (3)  and case (1) to case (4). In fact, $v\in I^{\mu}(N^*S_i)$ implies $v \in I^{\mu, l}(N^*S_{ij}, N^*S_i)$ for any $l$, when $S_i$ intersect $S_j$.   Therefore, it suffices to prove case (3) and (4). Indeed, one can work with the space 
\beq
\tilde \mci  \defeq  \sum_{1\leq i < j \leq J} I^{p, l}(N^*S_{ij}, N^*S_i)
\eeq
to show that this is an algebra and the terms \eqref{eqhigher} all belong  to $\tilde \mci.$ 

We also remark that it may happen that $S_{ij}  \cap S_{ik} \neq \emptyset$ for some $i < j < k$. This is the case when there is a line tangent to $D_i, D_j, D_k$. Hence the points in $S_{ij}, S_{kl}$ are not necessarily distinct.

We shall use the following local representation of paired Lagrangian distributions. We take $(x_1, x_2)$ to be local coordinates for $M$ near $0$ such that 
\beq
S_1 = \{x_1 = 0\}, \ \ S_2 = \{x_2 = 0\} \text{ and } S_{12} \defeq S_1\cap S_2 = \{x_1 = x_2 = 0\}.
\eeq
Let $\zeta = (\zeta_1, \zeta_2)$ be the dual coordinates of the cotangent space. Then for $u \in I^{p, l}(N^*S_{12}, N^*S_1)$ we have (see \cite[(1.5)]{GrU93})
\beqq\label{equexp}
u(x) = \int_{\mbr^2} e^{i x\cdot \zeta} a(x, \zeta) d\zeta
\eeqq
with $a$ belonging to the product type symbol space $S^{p, l-\ha}(\mbr^2\times \mbr\times \mbr)$ i.e.\ 
\beq
|\p_x^\gamma \p_{\zeta_2}^\beta \p_{\zeta_1}^\alpha a(x, \zeta)| \leq C_{\alpha\beta\gamma K} \langle \zeta \rangle^{p-|\alpha|} \langle \zeta_2\rangle^{l-\ha - |\beta|}, \ \ x\in K, 
\eeq
where $K\subset \mbr^2$ is any compact set and $C_{\alpha\beta K}$ is a positive constant.  
For $v \in I^{p, l}(N^*S_{12}, N^*S_2)$, we would have 
\beqq\label{eqvexp}
v(x) = \int_{\mbr^2} e^{i x\cdot \zeta} b(x, \zeta) d\zeta
\eeqq
with $b\in S^{p, l-\ha}(\mbr^2\times \mbr\times \mbr)$ and estimates
\beq
|\p_x^\gamma \p_{\zeta_2}^\beta \p_{\zeta_1}^\alpha b(x, \zeta)| \leq C_{\alpha\beta\gamma K} \langle \zeta \rangle^{p-|\alpha|} \langle \zeta_1\rangle^{l-\ha - |\beta|}, \ \ x\in K,
\eeq
for some positive constant $C_{\alpha\beta\gamma K}$. 
In fact, any $u\in I^{p, l}(N^*S_{ij}, N^*S_i)$ has a local representation as above, after conjugating by some elliptic FIO. See \cite{GU, DUV}. Thus it suffices to consider the pair $S_1, S_2$ above. 

We start with case (3). 
\begin{lemma}\label{lm2}
If $u, v \in I^{p, l}(N^*S_{12}, N^*S_1)$ and both $u, v$ are supported near $S_{12}$, then
\beq
uv \in I^{p, l}(N^*S_{12}, N^*S_1).
\eeq
\end{lemma}
\bpf
We use the local expression of $u(x)$ in \eqref{equexp} and let 
\beq
v(x) = \int_{\mbr^2} e^{i x\cdot \zeta} b(x, \zeta) d\zeta
\eeq
where $b \in S^{p, l-\ha}(\mbr^2\times \mbr\times \mbr)$ is a product-type symbol. Then 
\beq
\begin{split}
u(x)v(x) &=  \int_{\mbr^2 \times \mbr^2} e^{i x\cdot \zeta} a(x, \zeta) e^{i x\cdot \eta} b(x, \eta)  d\zeta d\eta = \int_{\mbr^2} e^{ix\zeta} c(x, \zeta)d\eta,
\end{split}
\eeq
where 
\beq
c(x, \zeta) = \int_{\mbr^2} a(x, \zeta - \eta) b(x, \eta) d\eta
\eeq
is a convolution. 
We consider the symbol estimates of $c(x, \zeta)$. We have
\beq
\begin{split}
&| \p_x^\gamma \p_{\zeta_2}^\beta \p_{\zeta_1}^\alpha c(x, \zeta)| \\
&\leq  C \int_{\mbr^2}  \langle \zeta_1-\eta_1, \zeta_2 -\eta_2 \rangle^{p-|\alpha|} \langle \zeta_2 - \eta_2 \rangle^{l-\ha-|\beta|} \cdot   \langle \eta  \rangle^{p} \langle \eta_2 \rangle^{l-\ha}     d\eta\\
& = C \langle \zeta \rangle^{p-|\alpha|} \langle  \zeta_2 \rangle^{l-\ha-|\beta|} \int_{\mbr^2} \dfrac{ \langle \zeta_1-\eta_1, \zeta_2 -\eta_2 \rangle^{p-|\alpha|}}{\langle \zeta \rangle^{p-|\alpha|}} \cdot \dfrac{\langle \zeta_2 - \eta_2 \rangle^{l-\ha-|\beta|}}{\langle  \zeta_2 \rangle^{l-\ha-|\beta|} } \cdot   \langle \eta  \rangle^{p} \langle \eta_2 \rangle^{l-\ha}     d\eta\\
&\leq C \langle \zeta \rangle^{p-|\alpha|} \langle  \zeta_2 \rangle^{l-\ha-|\beta|} \int_{\mbr^2} \langle \eta  \rangle^{p} \langle \eta_2 \rangle^{l-\ha}     d\eta
\end{split}
\eeq
We get the last inequality because the two fractions in the integrals are bounded for all $\zeta, \eta$. Finally, we see that since $p\leq -\frac 32$ and $l \leq -1$, the above integral converges. So we get that $c(x, \zeta) \in S^{p, l-\ha}(\mbr^2\times \mbr\times \mbr)$ and this concludes the proof. 
\epf

Using this lemma and \eqref{eqipl}, we immediately obtain the following 
\begin{cor}
  If $u \in I^{p, l}(N^*S_{12}, N^*S_1), v\in I^{\mu}(N^*S_1)$ and $v$ is supported near $S_{12}$, then
\beq
uv \in I^{p, l}(N^*S_{12}, N^*S_1).
\eeq
\end{cor}
\begin{cor}\label{lmprod}
If $u, v \in I^{\mu}(N^*S_1)$, then $uv \in  I^{\mu}(N^*S_1).$
\end{cor}

Next we consider case (4).
\begin{lemma}\label{lm3}
If $u \in I^{p, l}(N^*S_{12}, N^*S_1), v\in I^{\mu}(N^*S_{12}, N^*S_2)$ and both $u, v$ are supported near $S_{12}$, then
\beq
uv \in I^{p, l}(N^*S_{12}, N^*S_1) + I^{p, l}(N^*S_{12}, N^*S_2).
\eeq
\end{lemma}
\bpf
We take $u(x)$ as in \eqref{equexp} and let 
\beq
v(x) = \int_{\mbr^2} e^{i x\cdot \zeta} b(x, \zeta) d\zeta
\eeq
where $b$ is a product-type symbol satisfying
\beq
|\p_x^\gamma \p_{\zeta_1}^\beta \p_{\zeta_2}^\alpha b(x, \zeta)| \leq C_{\alpha\beta\gamma K} \langle \zeta \rangle^{p-|\alpha|} \langle \zeta_1\rangle^{l-\ha - |\beta|}.
\eeq
Then we have 
\beq
\begin{split}
u(x)v(x) &=  \int_{\mbr^2 \times \mbr^2} e^{i x\cdot \zeta} a(x, \zeta) e^{i x\cdot \eta} b(x, \eta)  d\zeta d\eta = \int_{\mbr^2} e^{ix\zeta} c(x, \zeta)d\eta,
\end{split}
\eeq
 where 
\beq
c(x, \zeta) = \int_{\mbr^2} a(x, \zeta - \eta) b(x, \eta) d\eta. 
\eeq

Now we let $\psi(t)\in C_0^\infty(\mbr)$ be a cut-off function so that $\psi(t) = 1$ for $|t|\leq \ha$ and $\psi(t) = 0$ for $|t|\geq 1$. Consider the symbol estimates for $\psi(\frac{\langle \zeta_2 \rangle}{\langle \zeta_1 \rangle}) c(x, \zeta)$ which is supported on $\langle \zeta_2\rangle \leq \langle \zeta_1 \rangle$. First we have
\beq
\begin{split}
|\psi(\frac{\langle \zeta_2 \rangle}{\langle \zeta_1 \rangle}) c(x, \zeta)| &\leq  \psi(\frac{\langle \zeta_2 \rangle}{\langle \zeta_1 \rangle}) \int_{\mbr^2} C \langle \zeta_1-\eta_1, \zeta_2 -\eta_2 \rangle^{p} \langle \zeta_2 - \eta_2 \rangle^{l-\ha} \cdot \langle  \eta  \rangle^{p} \langle \eta_1 \rangle^{l-\ha} d\eta \\
 & \leq  C \langle \zeta \rangle^{p} \langle \zeta_2 \rangle^{l-\ha} \int_{\mbr^2}\langle  \eta  \rangle^{p} \langle \eta_1 \rangle^{l-\ha} d\eta \leq C \langle \zeta \rangle^{p} \langle \zeta_2 \rangle^{l-\ha}.
 \end{split}
\eeq
Here we used the fact $l-\ha < 0$ to get $\langle \zeta_1\rangle^{l-\ha} \leq \langle \zeta_2 \rangle^{l-\ha}$ and also some consideration in the proof of Lemma \ref{lm2} to get the boundedness of the integral. Similarly, one can check the derivatives and conclude that $\psi(\frac{\langle \zeta_2 \rangle}{\langle \zeta_1 \rangle}) c(x, \zeta)  \in S^{p, l-\ha}(\mbr^2\times \mbr\times \mbr)$  and 
\beq
\int_{\mbr^2} e^{i x\cdot \zeta} \psi(\frac{\langle \zeta_2 \rangle}{\langle \zeta_1 \rangle}) c(x, \zeta) d \zeta \in I^{p, l}(N^*S_{12}, N^*S_1).
\eeq
By the same argument, we can show that $(1-\psi(\frac{\langle \zeta_2 \rangle}{\langle \zeta_1 \rangle}) c(x, \zeta)) \in S^{p, l-\ha}(\mbr^2\times \mbr\times \mbr)$ and 
\beq
\int_{\mbr^2} e^{i x\cdot \zeta} (1-\psi(\frac{\langle \zeta_2 \rangle}{\langle \zeta_1 \rangle})) c(x, \zeta) d \zeta \in I^{p, l}(N^*S_{12}, N^*S_2).
\eeq
This completes the proof. 
\epf

Of course, this lemma and \eqref{eqipl} imply that
\begin{cor} 
If $u \in I^{p, l}(N^*S_{12}, N^*S_1), v\in I^{\mu}(N^*S_2)$ and $v$ is supported near $S_{12}$, then
\beq
uv \in I^{p, l}(N^*S_{12}, N^*S_1) + I^{p, l}(N^*S_{12}, N^*S_2).
\eeq
\end{cor}
We also remark that Lemma \ref{lm3} and \eqref{eqipl}  imply \cite[Lemma 1.1]{GrU93} for those exponents we consider here. 
With all  these lemmas, we conclude that
\begin{prop}
The distribution space $\mci$ defined in \eqref{defs} is an algebra. The terms $[(\mcr \chi_D)^2]^k, k = 1, 2, \cdots,$ defined in \eqref{eqhigher} belong to $\mci.$
\end{prop}

 Finally, using Proposition \ref{propfbp}, we obtain the main result of this section.
\begin{theorem}\label{main1}
Consider $f_{MA}$ in $\mce'(\mbr^2)$  defined in \eqref{eqfct}. Suppose that the metal regions $D$ have smooth boundaries as described in Section \ref{pre} and satisfy assumptions (A1). We have that away from $\Sigma = \p D$, the streaking artifacts $f_{MA} \in  I^{-2}(N^*\mcl)$ with $\mcl$ defined in \eqref{defL}.
\end{theorem}
This results improves the main result Theorem 3.3 of \cite{Seo}. It gives a quantitative description of the distribution representing the streaking artifacts. Of course, the remarks (regarding the strictly convexity assumption) in the end of Section \ref{secnon} apply here as well.

\section{Metal regions with piecewise smooth boundaries}\label{secpiece}

As in Section \ref{pre}, we suppose that $D_j$ are simply connected domains of $\mbr^2$ with closed boundary $\Sigma_j$. In this section, we assume that $\Sigma_j$ are piecewise smooth curves, namely
\begin{enumerate}
\item $\Sigma_j = \gamma_j(t)$ for some $\gamma_j \in C^0([a, b]; \mbr^2)$ such that $\gamma_j(a) = \gamma_j(b)$;
\item  there exists a partition $a = t_0< t_1 < \cdots < t_{N_j} = b$ such that $\gamma_j^i \defeq \gamma_j|_{[t_{i-1}, t_{i}]}, i = 1, \cdots, N_j, N_j\geq 1$ are smooth.
\end{enumerate}
We point out that $a, b$ and the partition above depend on $\gamma_j$, though the dependency is not showing up in the notations. We denote the end points or corner points by $n_j^i \defeq \gamma_j(t_i), i = 1, \cdots, N_j$ and the collection of such points by $\mcn_j$. Finally, we let $\mcn = \bigcup_{j = 1}^J\mcn_j$. 

For piece-wise smooth boundaries, the singularities of the characteristic functions of $D_j$ near the corners could be rather complicated. In addition to (1) and (2), we require that 
\begin{enumerate}
\item[(3)]  $\gamma_j^i$ is either strictly convex or a line segment;
\item[(4)]  $\gamma_j^i$ can be smoothly extended across $t_{i-1}$ and $t_{i}$ to some smooth curves $\tilde \gamma_j^i$. Moreover, $\tilde \gamma_j^i$ intersect $\tilde \gamma_j^k, k = i-1, i+1$ transversally at the end points  $\gamma_j(t_i)$ and $\gamma_j(t_{i-1})$.
\end{enumerate}
Notice that the extensions $\tilde \gamma_j^i$ are not unique, and (3) allows us to treat polygon shaped regions. 

\begin{figure}[htbp]
\centering
\vspace{-0.7cm}
\includegraphics[scale=.8]{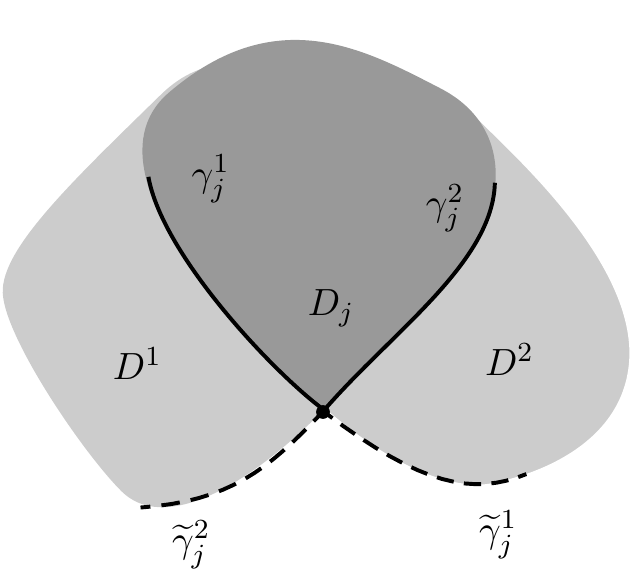}
\caption{$D_j$ is the region bounded by $\gamma_j^1, \gamma_j^2$. $D^1$ is the region above the curve $\tilde \gamma_j^2$ which extends $\gamma_j^2$. $D^2$ is the region above the curve $\tilde \gamma_j^1$ which extends $\gamma_j^1$.}
\label{figext}
\end{figure}

Under the above assumptions (1)-(4), we have
\begin{lemma}
The characteristic function of $D_j, j = 1, 2, \cdots, J,$ 
\beq
\chi_{D_j} \in \sum_{n_j^i\in \mcn_j} \big[ I^{-1, -\ha}(N^*n_j^i, N^*\gamma_j^i) + I^{-1, -\ha}(N^*n_j^i, N^*\gamma_j^{i+1})  \big].
\eeq
\end{lemma}
\bpf
It suffices to prove this locally near any corner points $n_j^i\in \mcn_j$. By (3) and (4) above, we can extend $D_j$ near $n_j^i$ to $D^1$ with boundary $\tilde \gamma_j^i$ and $D^2$ with boundary $\tilde \gamma_j^{i-1}$, see Figure \ref{figext}. Then near $n_j^i$, $\chi_{D_j} = \chi_{D_j^1}\cdot \chi_{D_j^2}$. Since $\chi_{D_j^k}, k = 1, 2$ are in $I^{-1}(N^*\tilde \gamma_j^i)$ and the boundaries intersect transversally, we get that 
\beq
\chi_{D_j} \in I^{-1, -\ha}(N^*n_j^i, N^*\gamma_j^i) + I^{-1, -\ha}(N^*n_j^i, N^*\gamma_j^{i-1}). 
\eeq
This finishes the proof.
\epf

Next consider $\mcr \chi_{D_j}, j = 1, \cdots J$. Since $\mcr$ is an elliptic Fourier integral operator, by the theory of paired Lagrangian distributions \cite{GU},  the composition $\mcr\chi_{D_j}$ are also paired Lagrangian distributions. To describe them, it suffices to find the intersecting Lagrangians after the symplectic transformation induced by $\mcr$. 

\begin{figure}[htbp]
\centering
\includegraphics[scale=.7]{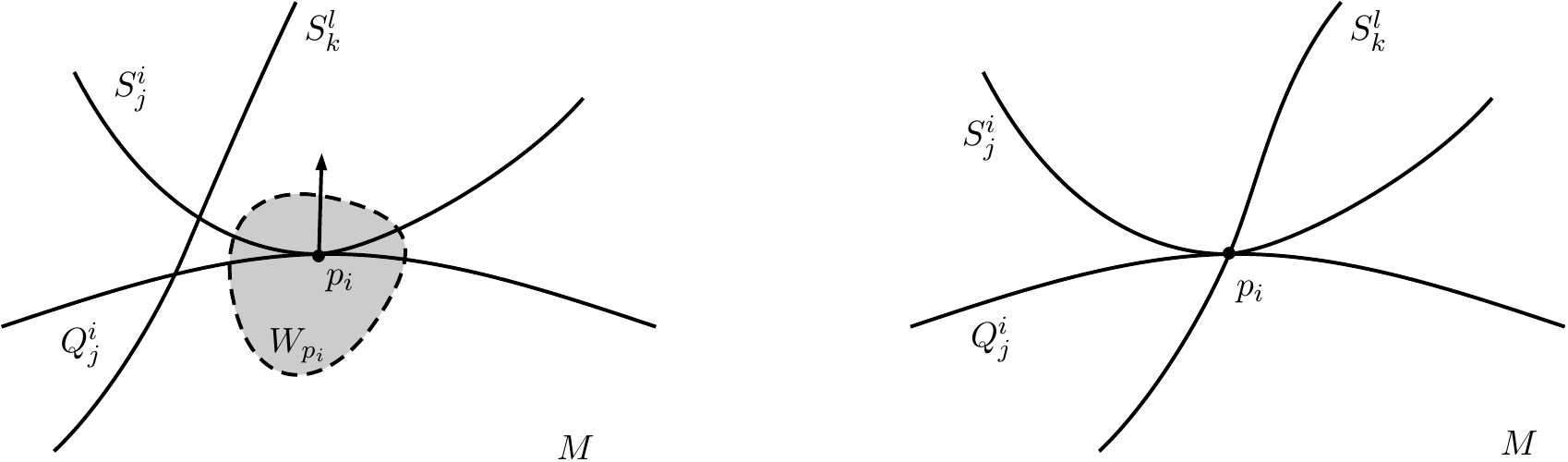}
\caption{Schematic view of possible intersections of curves in a sinogram. \textbf{Left:} Two curves $S_j^i$ and $Q_j^i$ intersects tangentially at $p_i$, and all the intersections with other $Q_k^l, S_k^l$ are outside of $W_{p_i}$. We introduce cut-offs supported in $W_{p_i}$.  \textbf{Right:} Another curve $S_k^l$ intersect  the two cruves at the point $p_i$. The point $p_i$ becomes a point singularity. }
\label{figtan}
\end{figure}

We recall the canonical relation $C$ of $\mcr$ defined in \eqref{eqcaC}. The composition 
\beq
 C\circ  N^*\gamma_j^i = \{(s, \theta; -\la, \la   x) \in T^*M\backslash 0 :   s = x\cdot \theta, \ \ x\in \gamma_j^i, \la \in \mbr\backslash 0, \theta \in N_x^*\gamma_j^i \cap \mbs^1\}.
\eeq
If $\gamma_j^i$ are strictly convex, the projection to $M$ is a one dimensional submanifold  of $M$ as seen in Section \ref{secnon}. For $j = 1, \cdots, J, i = 1, \cdots, N_j$, we define
\beq 
S_j^i \defeq \{(s, \theta) \in M : s = x\cdot \theta, \ \ \theta\in N_x^*\gamma_j^i \cap \mbs^1, x\in \gamma_j^i \},
\eeq 
then we have $C\circ  N^*\gamma_j^i =  N^*S_j^i$. If $\gamma_j^i$ is a line segment, then the projection of $C\circ  N^*\gamma_j^i$ to $M$ is $S_j^i$  which is a point. We still have $C\circ  N^*\gamma_j^i =  N^*S_j^i$.  Next, the composition 
\beq
 C\circ N^*n_j^i= \{(s, \theta; -\la, \la   x):   s =  x \cdot \theta, \ \ x = n_j^i, \la \in \mbr\backslash 0, \theta \in  \mbs^1\}.
\eeq
For $j = 1, \cdots, J, i = 1, \cdots, N_j$, we define
\beqq\label{eqqij}
Q_j^i \defeq \{(s, \theta)\in M : s = n_j^i \cdot \theta, \ \ \theta\in   \mbs^1 \}.
\eeqq 
We see that this is a one dimensional submanifold of $M$ and $C\circ N^*n_j^i = N^*Q_j^i$. The two Lagrangians $N^*S_{j}^i$ and $N^*Q_j^i$ intersect cleanly at a co-dimension one submanifold. Therefore, we proved 
\beqq\label{eqrdj}
\mcr\chi_{D_j} \in \sum_{i = 1}^{N_j}  \big[ I^{-\frac 32, -\ha}(N^*Q_j^i, N^*S_j^i) + I^{-\frac 32, -\ha}(N^*Q_j^i, N^*S_j^{i-1})\big], \ \ j = 1, 2, \cdots, J.
\eeqq
When $S_{j}^i$ is a point, $Q_j^i$ intersect $S_j^i$ transversally and we can analyze the singularities as before. However, when $S_j^i$ is a one dimensional submanifold, one can check that $S_j^i$ intersect $Q_j^i$ in a tangential way. This causes additional difficulties because the distributions in \eqref{eqrdj} do not fit in the framework we used in Section \ref{secfull}. Fortunately,  away from the intersection of $S_j^i$ and $Q_j^i$, we have a clear picture
\beq
\mcr\chi_{D_j} \in \sum_{i = 1}^{N_j}  \big[ I^{-2}(N^*Q_j^i)  + I^{-\frac 32}(N^*S_j^{i})\big]
\eeq
and this can be handled. More precisely, suppose $S_j^i$ is one dimensional and let $p_i \in S_j^i\cap Q_j^i$ be an intersection point. In fact, $p_i = (s_i, \theta_i)$ where $\theta_i$ is the unit normal vector at $n_j^i$ to $\gamma_j^i$ and $s_i = n_j^i\cdot \theta_i$.  We can find arbitrarily small open neighborhoods  $U_{p_i}, W_{p_i}$ of $p_i$ such that $\overline{U}_{p_i}\subset \subset W_{p_i}$. Now we introduce a smooth cut-off function $\psi_j^i$ so that $\psi_j^i = 1$ on $U_{p_i}$ and $\psi_j^i = 0$ outside of $W_{p_i}$. See the left of Figure \ref{figtan}. Let $\psi_j = \sum_i \psi_j^i$. Then we can write
\beqq\label{eqdec}
\begin{gathered}
\mcr\chi_{D_j} =  (1 - \psi_j) \mcr \chi_{D_j} + \psi_j \mcr\chi_{D_j} = \Psi_j^1 + \Psi_j^2 \text{ such that }\\
\Psi_j^1\defeq (1 - \psi_j) \mcr \chi_{D_j} \in \sum_{i=1}^{N_j} [I^{-2}(N^*Q_j^i)  + I^{-\frac 32}(N^*S_j^{i})] \text{ and }  \\
\Psi_j^2 \defeq \psi_j \mcr\chi_{D_j} \text{ is a paired Lagrangian distribution supported in $W_j \defeq \bigcup_{i = 1}^{N_j} W_{p_i}$}.
\end{gathered}
\eeqq
We shall analyze the nonlinear interactions of the terms $\Psi_j^1, j = 1, 2, \cdots, J$ as before. The nonlinear interactions of $\Psi_j^2$ are not so clear and we shall estimate the upper bound of their wave front sets.

\begin{lemma}\label{lml13}
We assume below that $j, k = 1, \cdots, J, i = 1, \cdots, N_j, l = 1, \cdots, N_l$, $(j, i)\neq (k, l)$,  and that the $S_j^i, S_k^l$ are one dimensional submanifolds.
\begin{enumerate}
\item If $S_j^i\cap S_k^l \neq \emptyset$, then they  intersect   transversally at a finite point set $S_{jk}^{il} \in M$. For each $p \in S_{jk}^{il}$,  there is a line $L_{p}$ tangent to both $\gamma_j^i$ and $\gamma_k^l$. The set of such tangent lines is denoted by $\mcl^1$. 
\item If $Q_j^i \cap Q_k^l \neq \emptyset$, then they intersect transversally at a finite point set $Q_{jk}^{il}$ and there is a (unique) line $L$ through $n_j^i$ and $n_k^l$. The set of such lines is denoted by $\mcl^2$.
\item If $S_j^i \cap Q_k^l \neq \emptyset$, then they intersect transversally at a finite point set $SQ_{jk}^{il}$. For each $p \in SQ_{jk}^{il}$, there is a line $L_p$ through $n_k^l$ and tangent to $\gamma_j^i$. The set of such lines is denoted by $\mcl^3.$
\end{enumerate}
\end{lemma}

\bpf
(1): This is the same as Lemma \ref{lmgeo}.

(2): In view of \eqref{eqqij}, we conclude that if $(s_0, \theta_0)\in Q_j^i\cap Q_k^l$, then $n_j^i$ and $n_k^l$ must lie on a straight line $x\cdot \theta_0 = s_0$. A similar argument as in Lemma \ref{lmgeo} shows that the intersection is transversal. 

(3): Suppose that $p  = (s, \theta) \in S_j^i\cap Q_k^l \neq \emptyset$. Then there exists $x_j \in \gamma_{j}^i$ and $x_k =n_k^l$ so that 
\beq
s = \theta\cdot x_j \text{ and } s = \theta\cdot x_k. 
\eeq
Therefore, $x_j, x_k$ lie on the straight line $L_p = \{x\in \mbr^2: s = \theta\cdot x\}$  tangent to $\gamma_j^i$ and intersects $n_k^l$. Again, a similar argument as in Lemma \ref{lmgeo} shows that the intersection is transversal. 
\epf

\begin{figure}[htbp]
\centering
\includegraphics[scale=.7]{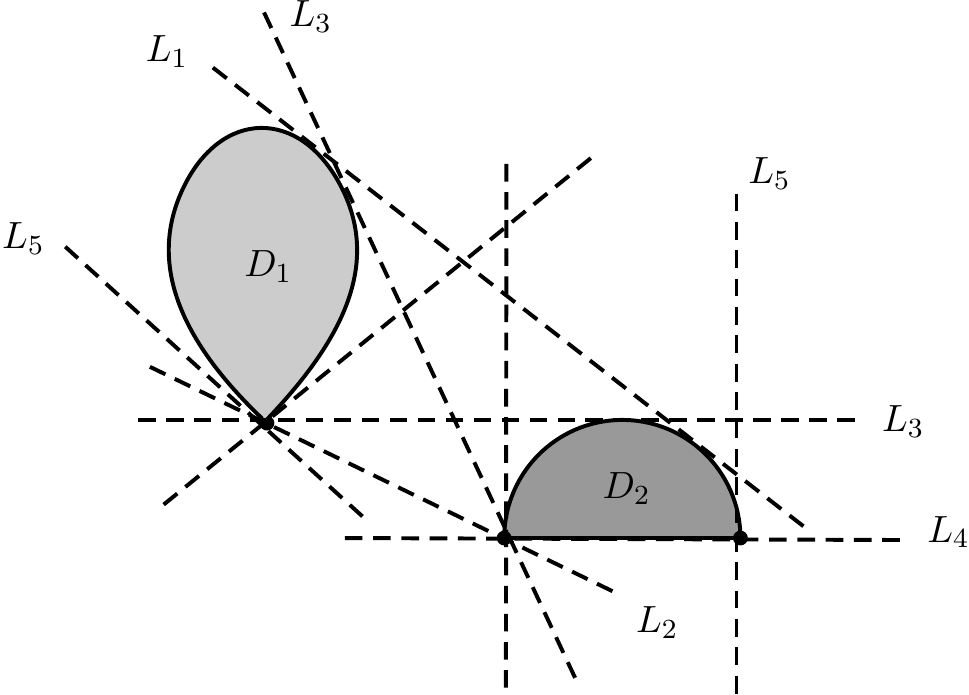}
\hspace{2em} \includegraphics[scale=.45]{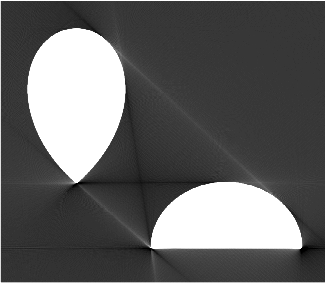}
\caption{Examples of streaking artifacts. \textbf{Left:} Theoretical results. The dashed lines represent the possible streaking artifacts. $L_1$ is tangent to two smooth curves of the two metallic objects $D_1, D_2$. $L_2$ is due to two corner points. $L_3$ is due to one corner point and one smooth curve. $L_4$ is due to the line segment of the boundaries and $L_5$ is due to the corner points. \textbf{Right:} Numerical results.}
\label{figstr}
\end{figure}

Now we are ready to discuss the nonlinear effects in $f_{MA}.$ We first write $\mcr(\chi_{D_j}) = \Psi_j^1 + \Psi_j^2, j = 1, 2, \cdots, J$ as defined in \eqref{eqdec} and we consider the nonlinear interactions of the $\Psi_j^1$ terms, away from the sets $W_j$ which contains intersection points $S_j^i\cap Q_j^i$.  Assume that $S_j^i, S_k^l, (j, i)\neq (k, l)$ are one dimensional submanifolds (not point sets). If $S_j^i, S_k^l$ intersect transversally at $p\in S_{jk}^{il}$, using \cite[Lemma 1.1]{GrU93}, we obtain that locally near $p$
\beq 
\mcr(\chi_{D_j})\mcr(\chi_{D_k}) \in I^{-\frac 32, -\frac 32 + \ha}(N^*p, N^*S_j^i) + I^{-\frac 32, -\frac 32+ \ha}(N^*p, N^*S_k^l).
\eeq 
If $S_{j}^i, Q_k^l$ intersect transversally at $p \in SQ_{jk}^{il}$, we obtain that 
\beq 
\mcr(\chi_{D_j})\mcr(\chi_{D_k}) \in I^{-\frac 32, -2 + \ha}(N^*p, N^*S_j^i) + I^{- 2, -\frac 32+ \ha}(N^*p, N^*Q_k^l).
\eeq 
If $Q_{j}^i, Q_k^l$ intersect transversally at $p \in Q_{jk}^{il}$, we obtain that 
\beq 
\mcr(\chi_{D_j})\mcr(\chi_{D_k}) \in I^{-2, -2 + \ha}(N^*p, N^*Q_j^i) + I^{- 2, -2+ \ha}(N^*p, N^*Q_k^l).
\eeq 
At this stage, one dimensional $S_j^i$ and $Q_k^l$ do not make a difference. So we let $\mcx$ be the set of the collection of such $S_j^i$ and $Q_k^l$. We relabel the element of $\mcx$ by $X_i, i = 1, \cdots, N$ and define a distributional space 
\beq
\mca \defeq \sum_{X_i, X_j\in \mcx} \sum_{p\in X_i\cap X_j} \big[ I^{-\frac 32, -1}(N^*p, N^*X_i) + I^{-\frac 32, -1}(N^*p, N^*X_j) \big].
\eeq
Then we have $\mcr(\chi_{D_j})\mcr(\chi_{D_k}) \in \mca$ away from the set $W_j, j = 1, 2, \cdots, J$. See the left of Figure \ref{figtan}. Of course, the orders above are not optimal. Actually, one sees that the artifacts due to different causes have different orders. For example, the artifacts due to the corner points are smoother than the ones in $\mcl^1.$ However, again our purpose is to describe the singularities in some distribution space which is an algebra. So we do not pursue a more precise statement. 

To complete the analysis, we consider the interactions in the set $W_j$, that is to consider $\Psi_j^2\cdot \Psi_k^2$ and $\Psi_j^1\cdot \Psi_k^2$. The latter term can be analyzed as before, because one can always shrink the neighborhoods $U_\bullet$ and $W_\bullet$ so that $S_j^i\cap Q_j^i$ does not intersect with $S_k^l$ or $Q_k^l$ in the support of $\Psi_j^1\cdot \Psi_k^2$, see Figure \ref{figtan}. It remains to consider $\Psi_j^2\cdot \Psi_k^2$. Recall that each $\Psi_j^2$ is a paired Lagrangian distribution, which is in $H^{s_0}(M)$ for some $s_0.$ The product $\Psi_j^2\cdot \Psi_k^2$ is a well-defined distribution. Provided that the supports of the functions are chosen to be sufficiently small, the wave front set is contained in the union 
\beq
\mcw \defeq \bigcup_{j = 1}^{J} \bigcup_{i = 1}^{N_j}[N^*S_j^i\cup N^*Q_j^i \cup N^*(S_j^i\cap Q_j^i)].
\eeq 
This is because away from $S_j^i\cap Q_j^i$, the distributions $\Psi_j^2$ are conormal and the product can be analyzed as in Section \ref{secfull}. 
We define a distributional space $\mcb \defeq \{u\in \mcd'(M): \WF(u)\subset \mcw\}$. Then what we just proved is $\mcr(\chi_{D_j})\mcr(\chi_{D_k}) \in \mca + \mcb$.

We are ready to prove the main result of this section. For $\mcl_i, i = 1, 2, 3$ defined in Lemma \ref{lml13}, we denote $\mcl' \defeq \bigcup_{i = 1}^3 \mcl_i. $ Next, we define $\mcl''$ to be the union of straight lines on $\mbr^2$ which are tangent to some $\gamma_j^i$ at the corner points. The main result below shows that the lines in $\mcl'$ and $\mcl''$ are all the possible streaking artifacts.

\begin{theorem}\label{main2}
Consider $f_{MA}$ in $\mce'(\mbr^2)$ defined in \eqref{eqfma} and suppose that the metal region  $D$ has piecewise smooth boundaries described in the beginning of this section.  We have that away from $\p D$, the streaking artifacts $f_{MA} \in I^{-2}(N^*\mcl')$ modulo a distribution whose wave front set is contained in $\mcl''$.
\end{theorem}
 
\bpf
By the algebraic property of $\mca$, we know that $[(\mcr \chi_D)^2]^k \in \mca + \mcb.$ Then
\beq
\mci^{-1}\big( [(\mcr \chi_D)^2]^k \big) \in \sum_{X_i, X_j\in \mcx} \sum_{p\in X_i\cap X_j} \big[ I^{-\frac 12, -1}(N^*p, N^*X_i) + I^{-\frac 12, -1}(N^*p, N^*X_j) \big] + \mcb.
\eeq
We can check that for all $p\in X_i\cap X_j$, $C^*\circ N^*p =  N^*L_p$, where $L_p$ is a line in $\mcl'$. The wave front set of the part in $\mce$ is contained in $\mcw$, and we see that $C^*\circ \mcw$ is contained in the union of $N^*\Sigma_j\cup N^*\mcn$ i.e.\ $\WF(\chi_D)$ and $\mcl''$. So we get 
\beq
\mcr^*\circ \mci^{-1}\big( [(\mcr \chi_D)^2]^k \big) \in \sum_{L_p\in \mcl'} \big[ I^{-1, -1}(N^*L_{p}, N^*\gamma_k^l) + I^{-1, -1}(N^*L_{p}, N^*\gamma_j^i) \big]
\eeq
modulo a distribution whose wave front sets is included in  $\WF(\chi_D)$ and $\mcl''$. This finishes the proof.  
\epf

\begin{figure}[htbp]
\centering
\includegraphics[scale=.8]{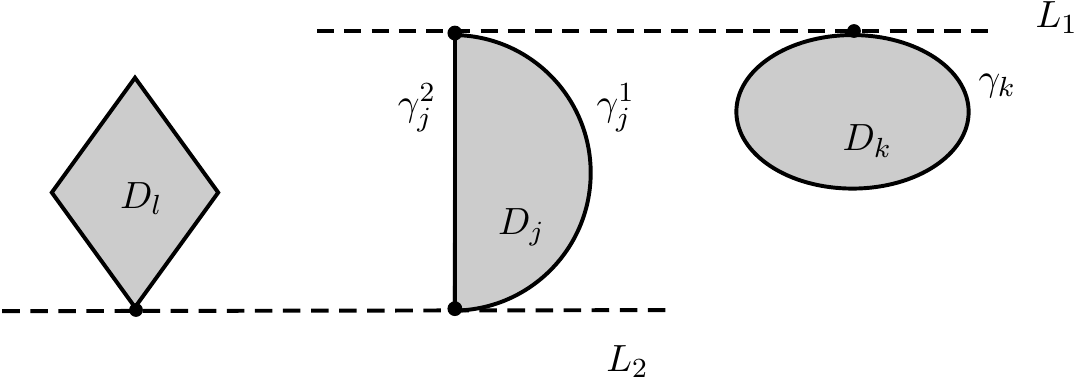}
\caption{ The line $L_1$ is tangent to $\gamma_k = \p D_k$ and passes through the corner point of $D_j$. Moreover, the line is tangent to $\gamma_j^1$. The sinogram is represented by the right of Figure \ref{figtan}. The line $L_2$ passes through two corner points and tangent to $\gamma_j^1$. These lines are also the possible streaking artifacts.}
\label{figa2}
\end{figure}

To conclude the section, we use two examples to illustrate the complexities of the artifacts for metal regions with piecewise smooth boundaries. 
First, we consider the possible interactions of singularities at $S_{j}^i\cap Q_{j}^i$ by other $S_k^l$ or $Q_k^l$. This type of geometry can be characterized easily. If $p \in S_{j}^i\cap Q_{j}^i\cap S_k^l$, there must be a line tangent to $\gamma_j^i$ and $\gamma_k^l$ and passing through the corner point on $\gamma_j^i$. If $p \in S_{j}^i\cap Q_{j}^i\cap Q_k^l$, then there is a line passing through a corner point on $D_k$ and tangent to $\gamma_j^i$, see Figure \ref{figa2}. In this case, one can see that the wave front set of $\mcr(\chi_{D_j})\mcr(\chi_{D_k})$ at $p$ is two dimensional. Therefore, the lines are expected to be present in the streaking artifacts. However, it is not clear whether they are conormal distributions.

\begin{figure}[htbp]
\centering
\hspace{-3em}
\begin{tabular}{cc}
\includegraphics[scale=.65]{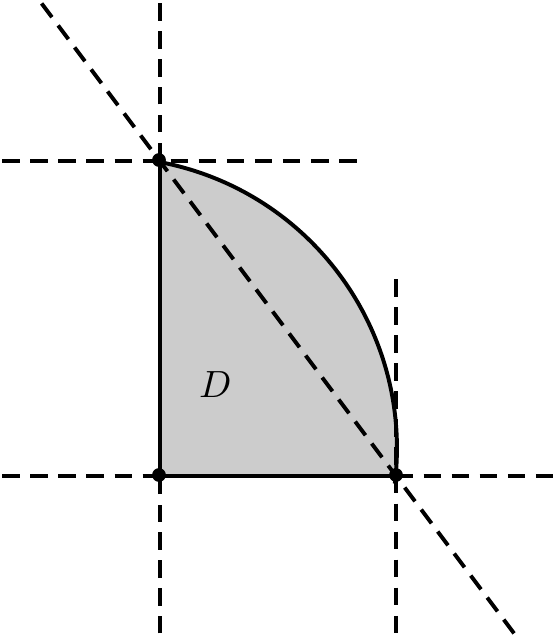}
& \includegraphics[scale=.4]{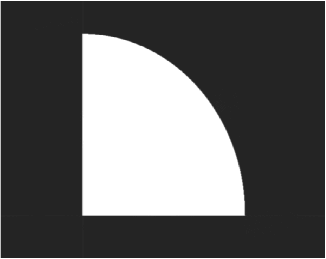}\\
\phantom{a} & $\mcr^*\mci^{-1}\mcr\chi_D$ \\
& \\
\includegraphics[scale=.4]{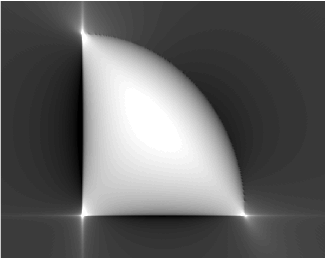} & \includegraphics[scale=.4]{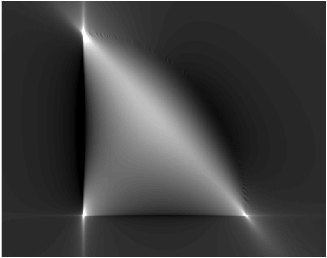}\\
 $\mcr^*\mci^{-1}[\mcr\chi_D]^2$ & $\mcr^*\mci^{-1}[\mcr\chi_D]^4$ 
\end{tabular}
\caption{Artifacts arising from the corners. Notice that the term $\mcr^*\mci^{-1}\mcr\chi_D$ only contains the conormal singularities of $\chi_D$ and doesn't carry any streaking artifacts, however they appear in the higher order terms as expected.}
\label{figeg}
\end{figure}

Next, we consider  a metal region consisting of a  simply connected region $D$ with boundary $\Sigma = \p D$. If $\Sigma$ is smooth, we learned from Section \ref{secnon} and \ref{secfull} that $\WF([\mcr\chi_D]^2)\subset \WF(\mcr\chi_D)$ and there is no streaking artifacts. However, if $\Sigma$ is piecewise smooth with line segments as shown in Figure \ref{figeg},  there could be streaking artifacts from the line segments and the corners. 
This explains the artifacts in Figure 3.3 of \cite{Seo}.

\section{Reduction of the artifacts}\label{secred}
Based on our quantitative analysis of the streaking artifacts, we study the reduction using appropriate filters before taking the X-ray transform. After these filters, the attenuation coefficients of the objects look smoother to the X-ray transform, see Figure \ref{fig:example1}. We will show that the non-linear effects become weaker in the sense that the order of the conormal distribution associated with the streaking artifacts are reduced.  Let $K$ be a pseudo-differential operator of order $k<0$ on $\mbr^2$. We consider 
\beq
K f_{E} = Kf_{E_0} + \alpha K(E-E_0)\chi_D
\eeq
as the modified attenuation coefficients. By the same mechanism, we have that the data
\beqq\label{new_data}
\tilde P(s,\theta) = \mcr Kf_{E_0} + \tilde P_{MA}, \ \ \tilde P_{MA} =  \sum_{k = 1}^N \frac{(-1)^k}{k} [\sum_{n = 1}^N \frac{(\alpha \eps)^{2n}}{(2n+1)!}(\mcr K \chi_D)^{2n}]^k.
\eeqq
We shall regard $\mcr K$ as X-ray transform with a pre-filter. The effect of $K$ is demonstrated in the theorem below. For simplicity, we only consider the setup of smooth boundaries. 

\begin{theorem}\label{main3}
Let $K$ be a pseudo-differential operator of order $k<0$ on $M$. Assume the setup of Theorem \ref{main1}. Then the reconstruction 
\beq
 \tilde f_{CT} \defeq K^{-1}\mcr^*\mci^{-1}\tilde P(s,\theta) = f_{E_0} + \tilde f_{MA},
\eeq
where away from $\p D$, the streaking artifacts $\tilde f_{MA} \in I^{-2 +  k}(N^*\mcl)$ with $\mcl$ defined in \eqref{defL}.
\end{theorem}
\bpf
Notice that $K\chi_D\in I^{\mu + k}(N^*\Sigma)$. We follow the proof of Theorem \ref{main1} line by line to complete the proof. This is possible because all the lemmas in Section \ref{secfull} work  for the exponent $\mu+k$ for $k< 0$.
\epf

\begin{figure}[!htbp]
\centering
\begin{tabular}{c}
\includegraphics[scale=.4]{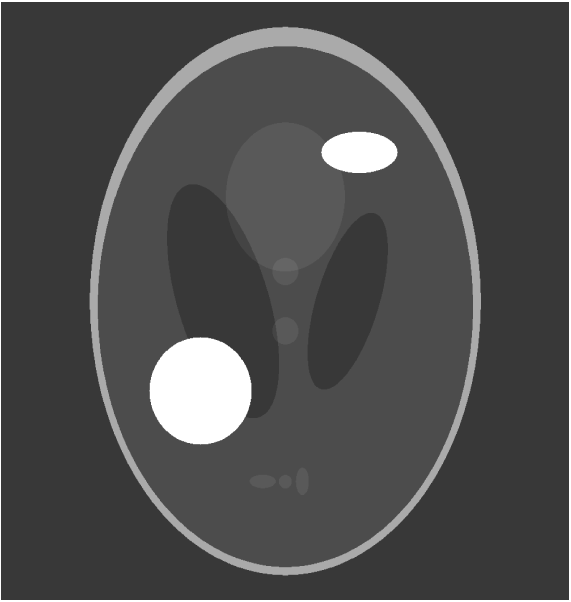}\quad\includegraphics[scale=.4]{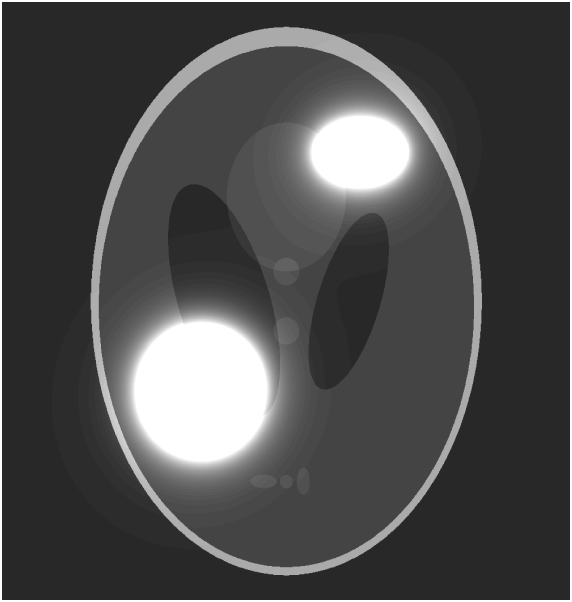}
\\
\includegraphics[scale=.4]{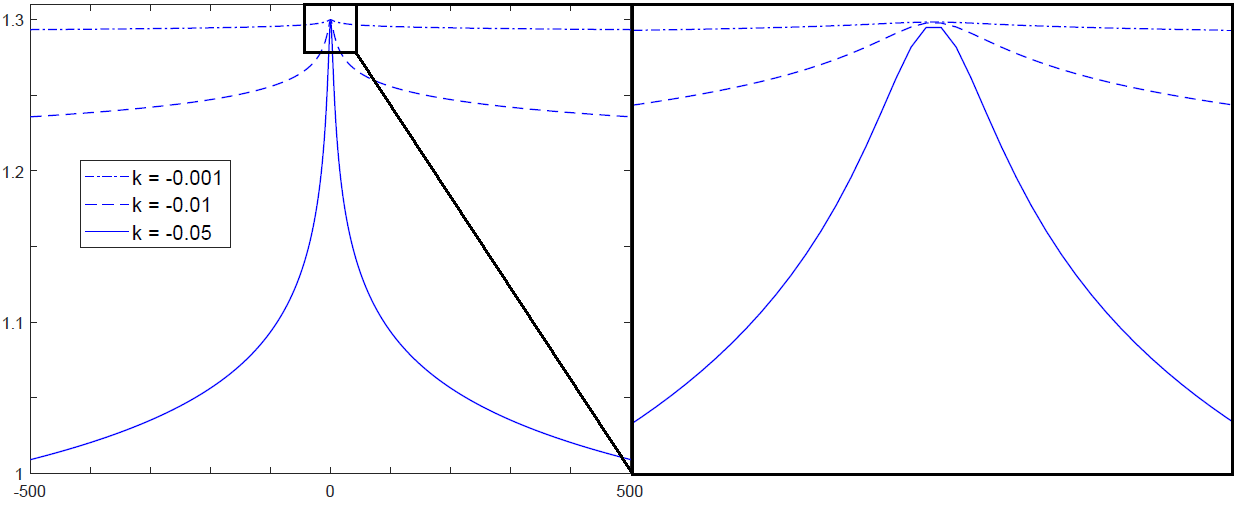}
\end{tabular}
\caption{(Top) Reference attenuation coefficient $f_{E_0}$ for numerical simulations in Figure \ref{fig:example1} and $Kf_{E_0}$ with $K$ as in \eqref{pre-filter} of order $k=-0.001$. The white regions correspond to the metal inclusions. (Bottom) Symbols of different orders for the filtering operator $K$ as in \eqref{pre-filter} with $\alpha=1.3$ and $k = -0.001$, $-0.01$ and $-0.05$.}
\label{fig:symbols}
\end{figure}
\begin{figure}[!htpb]
\centering
\includegraphics[scale=.4]{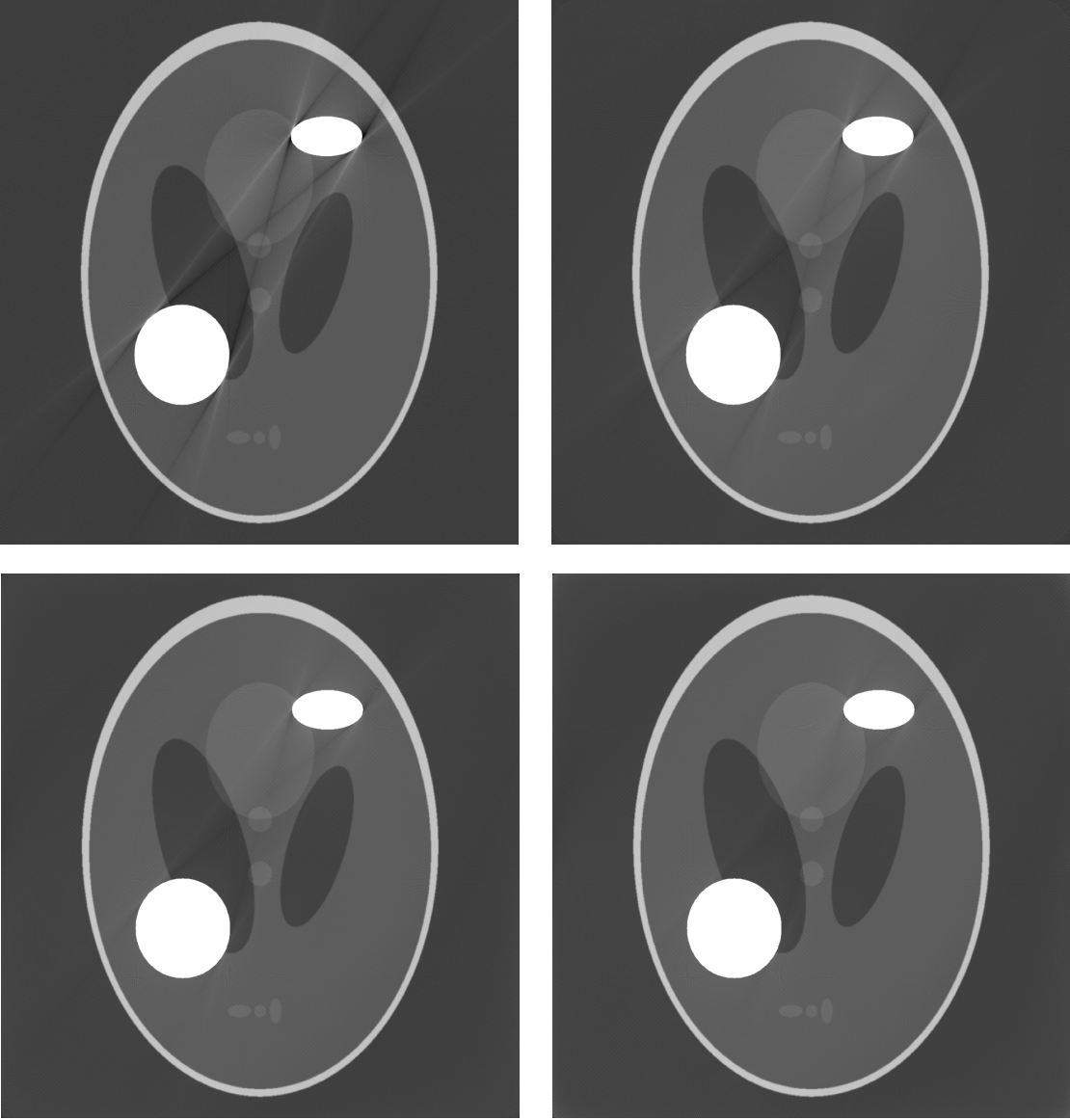}
\caption{(top left) FBP reconstruction $f_{CT}$. Proposed reconstruction $\tilde{f}_{CT}$ with filter $K$ given by the symbol \eqref{pre-filter} with $k=-0.001$ (top right), $-0.01$ (bottom left) and $-0.05$ (bottom right).}
\label{fig:example1}
\vspace{-0.5cm}
\end{figure}
 
\begin{figure}[!htbp]
\begin{center}
\includegraphics[scale=.4]{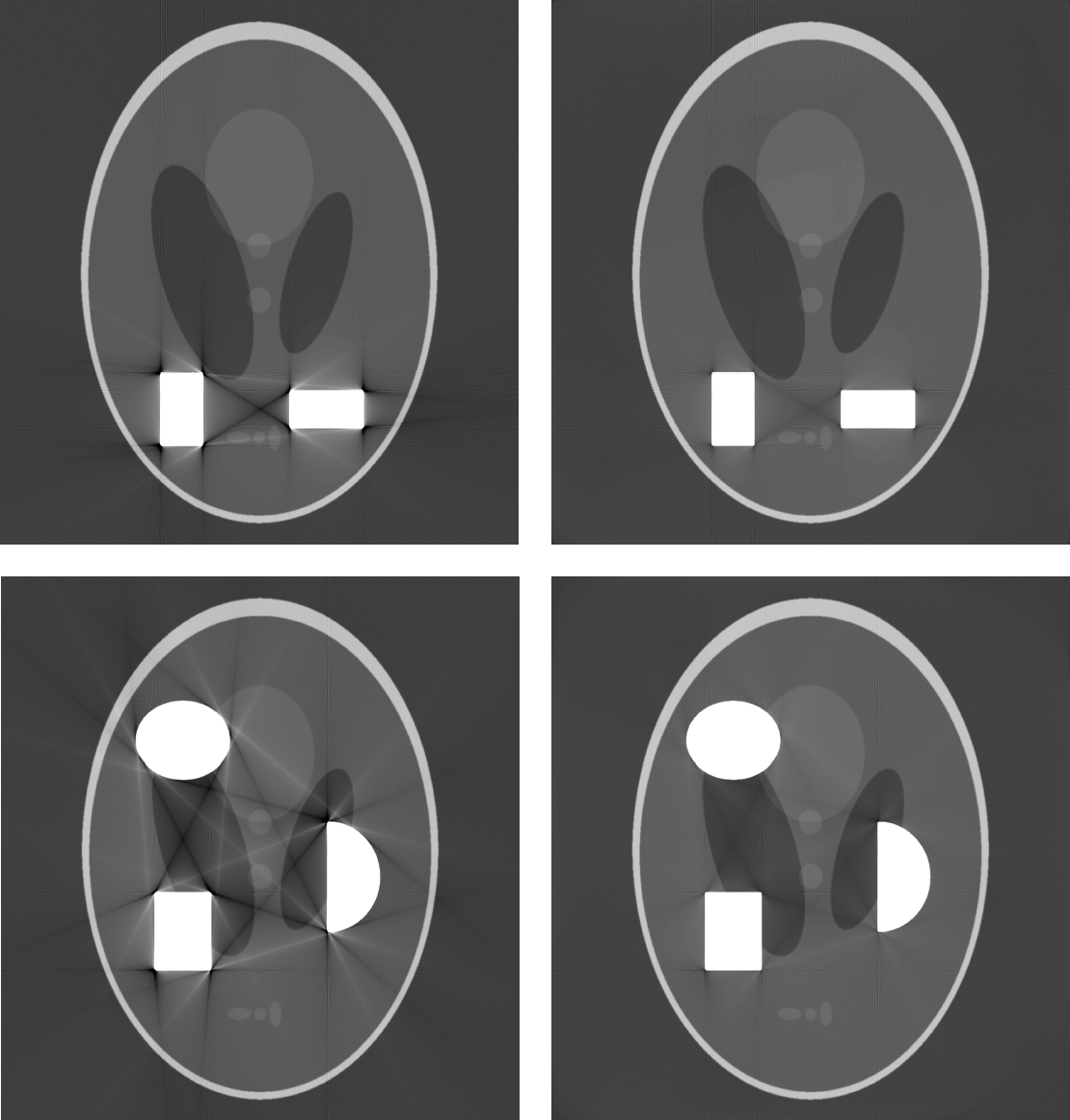} 
\end{center}
\caption{Reduction of the metal streaking artifacts by applying the filter $K$ of order $k=-0.01$. The left column corresponds to the CT image $f_{CT}$ and the right column is the enhanced reconstruction following Theorem \ref{main3}.}
\label{fig:example2}
\end{figure}

We remark that the reconstruction formula gives the same $f_{E_0}$ as without $K$. However, the streaking artifacts in $\tilde f_{MA}$ is smoother than $f_{MA}$ (the order of the conormal distribution is reduced by $k$). So it has the advantage that the reconstruction only reduces the streaking artifacts. It is also important to notice that in principle $K$ could be chosen as the identity times an small constant $\delta\in(0,1)$, this is, $K$ would act by homogeneously decreasing the attenuation coefficient, which then following \eqref{new_data} would imply in the reduction of the amplitude of the non-linear part and consequently of the streaking artifacts. However, we are interested in the reduction of artifacts through the interaction of  these filters with the order of the artifacts when considered as singularities. This is indeed the point of Theorem \ref{main3}. We will see in the examples below that even filters with small order of decay  reduce significantly the strength of the artifacts.
 
We conduct some numerical experiment to show the smoothing effects. We take operator $K$ to be a non-trivial pseudo-differential operator of order $k<0$ given by the symbol
\beqq\label{pre-filter}
p(x;\xi)=(\alpha+|\xi|^2)^{k/2},\quad (x;\xi)\in T^*\mbr^2 
\eeqq 
for some $\alpha>0$. In our examples we choose $\alpha=1.3$ so then $p(x,\xi)>1$ for all frequencies $|\xi|<500$. $K$ then corresponds to the elliptic operator $(\alpha-\Delta_x)^{k/2}$ that smooths out large frequencies of $f_E$. We let the computational grid to be of size $1000\times1000$ and the background attenuation coefficient of the soft tissue given by the Shepp-Logan phantom. Figure \ref{fig:symbols} shows the reference attenuation coefficient $f_{E_0}$ containing two metal regions with smooth boundary and the filtered reference coefficient $Kf_{E_0}$, as well as a cross section of $p(\cdot,\cdot)$ for different values of $k$. Notice that the smoothing effect of $K$ is stronger near the metal regions while the rest of the high frequencies remain almost intact.  
Figure \ref{fig:example1} contains a comparison of the FBP method and the reconstruction we propose for the different symbols in Figure \ref{fig:symbols}. 
In accordance with the theoretical analysis, the original (conormal) singularities of the absorption coefficient, corresponding to biological tissue and metal inclusions, remain intact while the streaking artifacts caused by the metal are reduced. This reduction depends on the shape of the symbol, being of course stronger when the symbol of the filter $K$ has a faster decay. In the cases of metal objects with piece-wise smooth boundaries, our method also reduces the artifacts generated by the corners as  can be seen in Figure \ref{fig:example2} where other examples of inclusion were considered.

 \appendix
 \section{Conormal and paired Lagrangian distributions}\label{secapp}
 For readers' convenience, we review the basics of these distributions and show some examples. The Lagrangian distributions can be found in H\"ormander \cite{Ho3, Ho4} and we refer the interested readers to \cite{DUV, FLU, GU, MU} for the detail about paired Lagrangian distributions. 
 
Let $\La$ be a closed conic Lagrangian submanifold of $T^*\mbr^n$. The space of Lagrangian distributions of order $\mu$, denoted by $I^\mu(\La)$, is defined as the set of all $u\in \mcd'(\mbr^n)$ such that 
\beqq\label{eqdefl}
L_1 L_2 \cdots L_N u \in  {}^\infty H^{\loc}_{-\mu - \frac n4}(\mbr^n)
\eeqq
for all $N$ and all properly supported first order pseudo-differential operators $L_j \in \Psi^1(\mbr^n)$ with principal symbols vanishing on $\La$, see \cite[Section 25.1]{Ho4}. Here, $ {}^\infty H^{\loc}_{-\mu - \frac n4}(\mbr^n)$ is the Besov space of the specified orders. Locally, such distributions can be expressed in terms of oscillatory integrals. Suppose that $\La$ is parametrized by a (homogeneous) phase function $\phi: \mbr^n\times \mbr^N\rightarrow \mbr$ over $U\subset \mbr^n$, namely $\La = \{(x, d_x\phi(x, \theta)) \in T^*\mbr^n: x \in U, d_\theta \phi(x, \theta) = 0\}$. Then $u\in I^\mu(\La)$ can be written as
\beq
u(x) = \int_{\mbr^N} e^{i \phi(x, \theta)} a(x, \theta) d\theta,
\eeq
where $a\in S^{\mu + \frac n4 - \frac N2}(\mbr^n\times \mbr^N)$ is a standard symbol satisfying 
\beq
|\p_x^\alpha \p^\beta_\theta a(x, \theta)| \leq C_{K \alpha \beta} \langle \theta \rangle^{\mu + \frac n4 - \frac N2 - |\beta|}, \ \ x\in K,
\eeq
for any compact set $K\subset U$, multi-indices $\alpha, \beta$ and some constant $C_{K, \alpha, \beta}.$ The wave front set of $u$ is contained in $\La$. The order $\mu$ is related to the strength of the distribution in terms of Besov regularities. One can also use Sobolev regularities $I^\mu(\La)\subset H^m(\mbr^n)$ for $m < -\mu - \frac n4.$  Let $K$ be a submanifold of $\mbr^n.$ The conormal bundle $N^*K$ is defined as 
\beq
 N^*K = \{(  x, \xi)\in T^*\mbr^n \backslash 0:   x\in K\ \ \xi|_{T_x K} = 0 \}.
\eeq
One can check that this is a conic Lagrangian submanifold of $T^*\mbr^n.$ The Lagrangian distributions $I^\mu(N^*K)$ are called conormal distributions. 

Conormal distributions naturally appear and are widely used in applications. Examples include the Heaviside functions and the Dirac delta function. As another example, we consider a homogeneous distribution $x_{1,+}^a, \re a > -1, x = (x_1, x_2, \cdots, x_n) \in \mbr^n$ defined as $x_{1,+}^a = x_1^a, x_1 > 0$ and $x_{1,+}^a = 0, x_1 \leq 0.$ The distribution has a conormal singularity at the hyper-surface $\{x\in \mbr^n: x_1 = 0\}.$ When $a = 0$, the distribution is just the Heaviside function (or the characteristic function) supported in $\{x_1\geq 0\}$. It is easy to see that 
\beq
x^a_{1, +} \cdot x_{1, +}^b = x_{1, +}^{a+b}, \ \ \re a > -1, \re b > -1.
\eeq
For $\re a$ or  $\re b > 0$, the distribution becomes smoother (e.g.\ in the sense of H\"older or Zygmund regularities) at $x_1 = 0$ after the multiplication, however, this won't be true if $a= 0 $ or $b= 0$. 
 
For two Lagrangian submanifolds $\La_0, \La_1$ of $T^*\mbr^n$ which intersect transversally at a codimension $k$ submanifold $\Omega$, the space of paired Lagrangian distributions of order $p, l$ is denoted by $I^{p, l}(\La_0, \La_1)$. Roughly speaking, these distributions have two types of singularities, the ones contained in $\La_0$ of order $p+l$, and the ones on $\La_1$ with order $p$.  In fact, for any $u\in I^{p, l}(\La_0, \La_1)$, we know that $u\in I^{p+l}(\La_0\backslash \Omega)$ and $u\in I^p(\La_1\backslash \Omega)$ as Lagrangian distributions. Also, $\bigcap_l I^{p, l}(\La_0, \La_1) = I^p(\La_1)$ and $\bigcap_p I^{p, l}(\La_0, \La_1) = C^\infty(\mbr^n).$ The paired Lagrangian distributions can be defined  using iterative applications of pseudo-differential operators as in \eqref{eqdefl}, see \cite{GrU}. Locally, they can be defined as oscillatory integrals. According to \cite[Proposition 2.1]{GU}, the clean intersecting Lagrangian pairs are locally equivalent under symplectmorphisms. We can find a local symplectmorphism $\chi$ so that the Lagrangian pair $(\La_0, \La_1)$ is transformed to the model pair $(N^*Y_2, N^*Y_1)$ where $Y_1 = \{x\in \mbr^n: x_1 = x_2 = \cdots = x_{d_1} = 0\} = \{x' = 0\}$ and $Y_2 = \{x\in \mbr^n: x_1 = x_2  = \cdots = x_{d_1 + d_2} = 0\} = \{x'' = 0\}$. (There are other model pairs, see \cite[Section 5]{DUV}.) Then we define $u \in I^{p, l}(N^*Y_2, N^*Y_1)$ by  
\beq
u(x) = \int_{\mbr^{d_1 + d_2}} e^{i(x'\cdot \xi' + x''\cdot \xi'')} a(x; \xi', \xi'')d\xi' d\xi''
\eeq
with $a(x; \xi', \xi'')$ belonging to the product-type symbol class 
\beq
S^{M, m}(\mbr^n; \mbr^{d_1}\backslash 0, \mbr^{d_2}) = \{a\in C^\infty: |\p_x^\gamma \p_{\xi''}^\beta \p_{\xi'}^\alpha a(x; \xi', \xi'')| \leq C_{K \alpha \beta \gamma} \langle \xi', \xi''\rangle^{M-|\alpha|} \langle \xi''\rangle^{m-|\beta|}\},
\eeq
where $M = p - \frac{d_1}{2} +  \frac n4$ and $m = l - \frac{d_2}{2}$. In general, we define distributions in $I^{p, l}(\La_0, \La_1)$ by $Fu, u \in I^{p, l}(N^*Y_2, N^*Y_1)$ with $F$  a zero order FIO associated with $\chi^{-1}$, see \cite[Definition 4.1]{GU}.

The paired Lagrangian distributions naturally appear in parametrix constructions for the (restricted) X-ray transforms \cite{FLU} and wave operators \cite{MU}. According to \cite[Lemma 1.1]{GrU93}, the product of two conormal distributions is a paired Lagrangian distribution if the two submanifolds intersect transversally. As an example, consider homogeneous distributions $x_{1,+}^a$ and $x_{2, +}^b$ in $\mbr^n, n\geq 2$. Then the product has a paired Lagrangian type singularity near $\{x_1 = x_2 = 0\}$.



\begin{thebibliography}{99}

\bibitem{BK} J. Barrett,  N. Keat. {\em Artifacts in CT: Recognition and avoidance.} Radio graphics : A Review Publication of the Radiological Society of North America, Inc, 24(6) (2004), 1679-91.

\bibitem{Seo0} J. K. Choi,  H. S. Park,  S. Wang,  Y. Wang,  J. K. Seo. {\em Inverse problem in quantitative susceptibility mapping.} SIAM Journal on Imaging Sciences, 7(3), (2014): 1669-1689.

\bibitem{FLU} D. Finch, I.-R. Lan, G. Uhlmann. {\em Microlocal analysis of the X-ray transform with sources on a curve.} Inside out: Inverse Problems and Applications 47 (2003): 193-218.

\bibitem{DUV} M. de Hoop, G. Uhlmann, A. Vasy. {\em Diffraction from conormal singularities.} Annales Scientifiques de l'\'Ecole Normale Sup\'erieure, 4e serie, t.48, (2015): 351-408.

\bibitem{GrU} A. Greenleaf,  G. Uhlmann. {\em Estimates for singular Radon transforms and pseudodifferential operators with singular symbols.} Journal of Functional Analysis 89.1 (1990): 202-232.

\bibitem{GrU93} A. Greenleaf, G. Uhlmann. {\em Recovering singularities of a potential from singularities of scattering data.} Communications in Mathematical Physics 157.3 (1993): 549-572.

\bibitem{GU} V. Guillemin, G. Uhlmann. {\em Oscillatory integrals with singular symbols.} Duke Math. J 48.1 (1981): 251-267.

\bibitem{Ho3} L. H\"ormander. {\em The analysis of linear partial differential operators III: Pseudo-differential operators.} Classics in Mathematics (2007).

\bibitem{Ho4} L. H\"ormander. {\em The analysis of linear partial differential operators IV: Fourier integral operators.} Classics in Mathematics (2009).

\bibitem{MU} R. Melrose, G. Uhlmann. {\em Lagrangian intersection and the Cauchy problem.} Communications on Pure and Applied Mathematics 32.4 (1979): 483-519.

\bibitem{Na} F. Natterer. {\em The mathematics of computerized tomography.} B. G. Teubner, Stuttgart; Wiley, Chichester, 1986.

\bibitem{PUW} B. Palacios,  G. Uhlmann, Y. Wang. {\em Reducing streaking artifacts in quantitative susceptibility mapping.}  SIAM Journal on Imaging Sciences Vol.10, No.4,  pp.1921-1934. (2017).

\bibitem{Seo} H. S. Park, J. K. Choi,  J. K. Seo. {\em Characterization of metal artifacts in X-ray computed tomography.} Communications on Pure and Applied Mathematics (2017): published on-line.

\bibitem{Seo1} H. S. Park, D. Hwang, J. K. Seo. {\em Metal artifact reduction for polychromatic X-ray CT based on a beam-hardening corrector.} IEEE Transactions on Medical Imaging 35.2 (2016): 480-487.

\bibitem{ZDL} H. Zhang, B. Dong, B. Liu. {\em A re-weighted joint spatial-Radon domain CT image reconstruction model for metal artifact reduction.} Preprint, 2017. 

\end{thebibliography}
\end{document}